\documentclass[12pt]{article}
\usepackage{epsf,amsmath, amsfonts}
\usepackage{epsfig,mathtools}
\usepackage{color,graphicx}% Include figure files
\usepackage{dcolumn}% Align table columns on decimal point
\usepackage{bm}% bold math
\usepackage{epsfig,latexsym}
\usepackage{lineno,newtxtext,nicefrac,bigints}
\usepackage{dashrule,hyperref}
\usepackage{mathptmx}      % use Times fonts if available on your TeX system

\usepackage[skip=0pt]{caption}
\usepackage[skip=0pt]{subcaption}
\newcommand{\bseq}{\begin{subequations}}
\newcommand{\eseq}{\end{subequations}}
\newcommand{\beq}{\begin{equation}}
\newcommand{\eeq}{\end{equation}}
\newcommand{\bef}{\begin{figure}}
\newcommand{\eef}{\end{figure}}

\hypersetup{
    colorlinks=true,
    urlcolor=red,
}

\let\OLDthebibliography\thebibliography
\renewcommand\thebibliography[1]{
  \OLDthebibliography{#1}
  \setlength{\parskip}{0pt}
  \setlength{\itemsep}{1pt plus 0.3ex}
}

\begin{document}
%\linenumbers
%\date{}

\title{\large \bf Implicit-explicit-compact methods for advection diffusion reaction equations} 
\author{ {\bf S. Singh$\dagger$} and {\bf S. Sircar$\dagger, \ddagger$}
\\{$\dagger$ \small Department of Mathematics, IIIT Delhi, India 110020} \\{\small $\ddagger$ Corresponding Author (email: sarthok@iiitd.ac.in)}}
\maketitle

\begin{abstract}
We provide a preliminary comparison of the dispersion properties, specifically the time-amplification factor, the scaled group velocity and the error in the phase speed of four spatiotemporal discretization schemes utilized for solving the one-dimensional (1D) linear advection diffusion reaction (ADR) equation: (a) An explicit (RK$_2$) temporal integration combined with the Optimal Upwind Compact Scheme (or OUCS3, J. Comp. Phys., 192, pg. 677-694 (2003)) and the central difference scheme (CD$_2$) for second order spatial discretization, (b) a fully implicit mid-point rule for time integration coupled with the OUCS3 and the Lele's compact scheme (J. Comp. Phys., 103, pg. 16-42 (1992)) for first and second order spatial discretization, respectively, (c) An implicit (mid-point rule)-explicit (RK$_2$) or IMEX time integration blended with OUCS3 and Lele's compact scheme (where the IMEX time integration follows the same ideology as introduced by Ascher et al., SIAM J. Numer. Anal., 32(3), pg. 797-823 (1995)), and (d) the IMEX (mid-point/RK$_2$) time integration melded with the New Combined Compact Difference scheme (or NCCD scheme, J. Comp. Phys., 228, pg. 6150-6168 (2009)). Analysis reveal the superior resolution features of the IMEX-NCCD scheme including an enhanced region of neutral stability (a region where the amplification factor is close to one), a diminished region of spurious propagation characteristics (or a region of negative group velocity) and a smaller region of nonzero phase speed error. The dispersion error of these numerical schemes through the role of $q-$waves is further investigated using the novel error propagation equation for the 1D linear ADR equation. Again, the in silico experiments divulge excellent Dispersion Relation Preservation (DRP) properties of the IMEX-NCCD scheme including minimal dissipation via implicit filtering and negligible unphysical oscillations (or Gibbs' phenomena) on coarser grids. The spectral resolution of the IMEX-NCCD scheme further is benchmarked by solving the classical two-dimensional (2D), Patlak-Keller-Segel (PKS) model. Numerical results reveal that the spiky structure of the solution is oscillation free and, when compared with the time-explicit OUCS3-CD2 method, the solution is better resolved by the IMEX-NCCD method.
\end{abstract}


%%%%%%%%%%%%%%%%%%%%%%%%%%%%%%%%%%%%%%%%%%%%%%%%
\section{Introduction} \label{sec:intro}
Models describing multiphysics and multiscale processes are ubiquitous in numerical simulations. Classical applications include those arising in aeronautics, meterology, biology, material and environmental sciences which are modeled by Navier-Stokes~\cite{TKS2013, Sircar2010}, reaction-diffusion~\cite{Sircar2016EPJ, Sircar2016DCDS} or ADR equations~\cite{Ascher1995, Sircar2015}. The individual physics or scale components have very different properties that are reflected in their discretization; for example, for ADR systems, the discrete advection has a relatively slow (or `nonstiff') dynamics while the diffusion and chemistry are fast (or `stiff') evolving processes~\cite{Hundsdorfer2003,Sircar2016JOMB}. The discretization in time of slow processes is advantageous with an explicit method since the nonstiff terms are (often) nonlinear and computationally more expensive to evaluate whereas implicit methods are more appropriate for stiff processes because of their favorable stability properties, allowing one to select longer time steps~\cite{Whitaker2013}. IMEX integrators have been proposed as an attractive alternative (compared with fully explicit or fully implicit time integration methods) where one combines the implicit (explicit) integration for the fast (slow) scale~\cite{Crouzeix1980}. IMEX linear multistep methods are investigated in~\cite{Hundsdorfer2007}, while IMEX Runge-Kutta schemes have been classified in~\cite{Ascher1997}. The literature on higher order IMEX methods are limited since these are difficult to construct due to large number of order constraints, increasing stability restrictions with increasing order of accuracy~\cite{Boscarino2013} as well as the malaise of the factitious `computational mode' that all multi-time level discretization schemes, with three or more time levels suffer~\cite{TKS2013}.

Similarly, the spatial discretization schemes must resolve all scales present in the flow for accurate numerical solution~\cite{Lighthill1993}. Additionally, the numerical propagation properties of every individually resolved scale (i.~e., the numerical group velocity and the numerical phase speed) must be the same as the physical propagation properties of the respective scale, so-called the DRP property~\cite{TKS2007}. The Fourier spectral methods are known for their ability to provide numerical accuracy and dispersion error-free results with all scales resolved and aliasing problems avoided~\cite{Gottlieb1977}. But the application of such methods are largely limited to simple geometries and boundary conditions with uniformly spaced nodes~\cite{Canuto1987}, with some isolated cases attempted in complex domians~\cite{Ditkowski2001,Pind2019}. Alternative to spectral method are the higher order explicit discretization methods~\cite{Tam1992} and methods based on Pad\'e approximation~\cite{Lele1992}, synonymous with the compact schemes. All compact schemes for evaluating the $n$th derivative of an unknown vector $\{ u \}$ require solving the set of linear algebraic equations given by $[A] \{u^{(n)}\} = [B]\{u\}$ or
\beq
u^{(n)} = [A]^{-1}[B]\{u\} = [D]\{u\}.
\label{eqn:CompStencil}
\eeq
Compact schemes have excellent resolution property~\cite{TKS2003}, even when their formal order of accuracy is lower than the traditional explicit discretization methods~\cite{TKS2006}. Although compact schemes have fewer points in the stencil, their implicit nature gives a global approximation of the derivatives~\cite{TKS2013}. Other relevant numerical properties of compact schemes have been analyzed earlier, including the DRP which is essential in solving time-dependent problems accurately ~\cite{TKS2009P1}, the de-aliasing or the unphysical pile-up in the energy spectrum for high wavenumbers~\cite{TKS2009P2} and higher spectral resolution on coarser meshes~\cite{TKS2010}. However, in this article we examine the dispersion properties of the combined, spatiotemporal IMEX-Compact schemes using the error propagation equation, described next.

Traditional stability analysis attributed to von Neumann for linear problems assume error to follow the same dynamics as given by the governing difference equation for the signal~\cite{Charney1950}. Recent studies have identified mechanisms for solution breakdown of linear equation to be related to the numerical properties of discrete computing, i.~e., the dissipation, dispersion and phase errors, and that the governing error propagation equation is different from the governing dynamics of the signal~\cite{TKS2007,TKS2014}. The numerical schemes for non-periodic problems discussed in \S \ref{sec:spectral}, require a similar analysis and for this purpose we have selected the 1D linear ADR equation,
\beq
\frac{\partial u}{\partial t} + c \frac{\partial u}{\partial x} = \nu \frac{\partial^2 u}{\partial x^2} + \lambda u, \quad c, \nu > 0.
\label{eqn:ADR}
\eeq
where $c, \nu, \lambda$ are constant advection speed, diffusivity and growth rate, respectively. The technique for analyzing discrete computational methods using Fourier-Laplace transforms is invoked in this article~\cite{TKS2007}, where the unknown is expressed in the wave number ($k$) plane by, $u(x, t) = \int U(k, t) e^{{\it} kx} dk$, such that an amplification factor can be introduced as $G(k) = \nicefrac{U(k, t+\triangle t)}{U(k,t)}$. For direct numerical simulations (DNS) one must have a neutrally stable method ($|G|=1$). If we define the numerical solution of equation~\eqref{eqn:ADR} by $u^N$ and the numerical error as $e = u - u^N$, then the correct error propagation, partial differential equation (PDE) is given by (refer \S \ref{sec:appendA} for the detailed derivation),
\begin{align}
&\frac{\partial e}{\partial t} + c \frac{\partial e}{\partial x} -\nu \frac{\partial^2 e}{\partial x^2} - \lambda e =\nonumber \\
&\bigintsss \!\!\!\!A_0(k) ({\it i} k c + \nu k^2 \!\!-\!\! \lambda) \!\!\!\left[ 1 \!\!+\!\! \frac{h}{c \triangle t}\!\!\left( \frac{\ln |G_{\text{num}}| \!+\! {\it i} \tan^{-1}\!\!\left[\frac{(G_{\text{num}})_{\text{imag}}}{(G_{\text{num}})_{\text{real}}}\right]}{\frac{\nu (k h)^2 - \lambda h^2}{c h} + {\it i}k h} \right) \!\!\right] \!\![G_{\text{num}}]^{t/\triangle t} e^{{\it i}k x} dk,
\label{eqn:ErrorADR}
\end{align}
where $A_0(k)$ denotes the Fourier amplitude of the initial condition, $G_{\text{num}}$ is the amplification factor of the numerical discretization method and $\triangle t, h$ are the temporal and the spatial step-size, respectively. One notes that equation~\eqref{eqn:ErrorADR} clubs error based on generic stability (through the numerical amplification factor) and the DRP properties of the method of discretization (e.~g. the scaled group velocity and the absolute error in the phase speed), unlike the one obtained using the modified equation approach~\cite{Gustafsson1972}, where the truncation error is accounted by collating and representing the discretized terms in the difference equation by their equivalent differential forms.

The aim of this article is to investigate the DRP properties of four, high accuracy (two-time level) spatiotemporal discretization using the state of the art error propagation equation. The two-time stepping methods are selected since they are devoid of the hazards of the contrived computational modes~\cite{TKS2013} and the complication related to the linear resonance stability~\cite{Stern2009}. \S \ref{sec:spectral} details the spectral analysis of these four numerical methods including the Explicit-OUCS3-CD$_2$ scheme (\S \ref{subsec:scheme1}), the Implicit-OUCS3-Lele scheme (\S \ref{subsec:scheme2}), the IMEX-OUCS3-Lele scheme (\S \ref{subsec:scheme3}) and the IMEX-NCCD scheme (\S \ref{subsec:scheme4}). \S \ref{sec:ED} examines the DRP properties of these four discretization techniques using the error propagation equation~\eqref{eqn:ErrorADR}. The spectral features of these four numerical methods (and in particular the IMEX-NCCD scheme) is benchmarked via the solution of the nonlinear, parabolic-elliptic PKS chemotaxis model in \S \ref{sec:PKS}. \S \ref{sec:conclusion} concludes with a brief discussion of the implication of these results. %as well as the focus of our future direction.

%%%%%%%%%%%%%%%%%%%%%%%%%%%%%%%%%%%%%%%%%%%%%%%%
\section{Spectral analysis for linear 1D ADR equation} \label{sec:spectral}
The basis of the spectral analysis for the numerical methods formulated for equation~\eqref{eqn:ADR}, is to transform the information from the physical to the spectral plane using its Fourier-Laplace transform. The spectral resolution of the discretization methods discussed below (i.~e., the absolute value of the ratio of the numerical to the physical amplification factor, $G=|\nicefrac{G_{\text{num}}}{G_{\text{exact}}}|$, the scaled group velocity, $V_g=\nicefrac{V_{g,\text{num}}}{V_{g,\text{phy}}}$, and the absolute value of the error in the phase speed, $c^{\text{err}}_{\text{phase}}=| 1 - \nicefrac{c_{\text{num}}}{c_{\text{exact}}} |$, refer \S \ref{sec:appendA}) is investigated for all wavenumbers within the range $k h \in [0, \pi]$ where the upper bound of this range is determined by the Nyquist criterion~\cite{Press1989}. $h$ is the grid spacing. The non-dimensional parameters introduced in the analysis are the Courant-Friedrichs-Lewy (CFL) number, $N_c = \frac{c \triangle t}{h}$, the Peclet number, $Pe = \frac{\nu \triangle t}{h^2}$, and the Damk\"ohler number, $Da = \lambda \triangle t$. For illustration purpose, $Pe=0.01$ and $Da=-0.01$ is fixed, the domain is discretized using $N=1001$ equidistant points and the plots in \S \ref{subsec:scheme1}-\ref{subsec:scheme4} are shown for an interior node, $m=500$.

%%%%%%%%%%%%%%%%%%%%%%%%%%%%%%%%%%%%%%%%%%%%%%%%
\subsection{Explicit-OUCS3-CD$_2$} \label{subsec:scheme1}
In the first numerical method, the temporal discretization of equation~\eqref{eqn:ADR} is achieved via an explicit two-stage Runge Kutta/Heun's method (also known as the explicit trapezoidal rule)~\cite{Ascher1995}. The advection term is spatially discretized using the OUCS3~\cite{TKS2003} and the diffusion component is discretized by the central CD$_2$ scheme.

Assuming a domain with $(N+1)$ equidistant points with spacing $h$ where a function $u$ is defined, the second order accurate OUCS3 estimates the first spatial derivative of the solution at the $j$th node, $u'_j$~\cite{TKS2013}. The inner stencil of the scheme is given by
\beq
p_{j-1} u'_{j-1} + u'_j + p_{j+1} u'_{j+1} = \frac{1}{h} \sum^{r=2}_{r=-2} q_r u_{j+r},
\label{appendB:interior}
\eeq
where $p_{j \pm 1} = D \pm \frac{\eta}{60}$, $q_{\pm 2} = \pm \frac{F}{4} + \frac{\eta}{300}$, $q_{\pm 1} = \pm \frac{F}{2} + \frac{\eta}{30}$, $q_0 = -\frac{11\eta}{50}$ with the coefficients $D = 0.3793894912, F=1.57557379, E = 0.183205192$ and $\eta = -2$. For an improved numerical stability, the one-sided boundary stencil ($j = 1, N+1$) and the near boundary stencil ($j = 2, N$) are proposed as follows,
\bseq
\begin{alignat}{4}
&u'_1 = \frac{1}{h}(-1.5 u_1 + 2 u_2 - 0.5 u_3), \label{appendB:up1} \\
&u'_{N+1} = \frac{1}{h}(1.5 u_{N+1} - 2 u_N + 0.5 u_{N-1}), \label{appendB:upN+1} \\
&u'_2 = \frac{1}{h} \left[ \left( \frac{2\beta_2}{3} - \frac{1}{3} \right) u_1 - \left( \frac{8\beta_2}{3} + \frac{1}{2} \right) u_2 + (4\beta_2 + 1) u_3 - \left( \frac{8\beta_2}{3} + \frac{1}{6} \right) u_4 + \frac{2\beta_2}{3} u_5 \right], \label{appendB:up2} \\
&u'_N = -\frac{1}{h} \left[ \left( \frac{2\beta_N}{3} - \frac{1}{3} \right) u_{N+1} - \left( \frac{8\beta_N}{3} + \frac{1}{2} \right) u_N + (4\beta_N + 1) u_{N-1} - \left( \frac{8\beta_N}{3} + \frac{1}{6} \right) u_{N-2} + \frac{2\beta_N}{3} u_{N-3} \right], \label{appendB:upN}
\end{alignat}
\label{appendB:boundary}
\eseq
where $\beta_2 = -0.025$ and $\beta_N = 0.09$. The point $j = 1$ in all of the compact schemes discussed in this and the next sections is of little concern since the governing PDF is not discretized at this node, only the boundary condition is specified there. In OUCS3, the point $j = 2$ shows instability across all wavenumbers and therefore this calculated derivative is eventually replaced by the CD$_2$ scheme locally. The solution of equation~\eqref{eqn:ADR} at the $j$th node ($u^*_j, u^{n+1}_j$) at each stage of the temporal discretization as a function of the solution at previous time-step, $u^n_j$, is outlined as follows,
\bseq
\begin{alignat}{4}
&u^{*}_j=u^{n}_j - N_c \sum_{r=1}^{N+1}(D)_{jr}u_{jr}^n + Pe\bigg( u_{j-1}^n-2 u_j^n+u_{j+1}^n \bigg) + (Da) u^n_j, \label{appendB:u*}
\\
&u^{n+1}_j=u^{n}_j - \frac{Nc}{2} \sum_{r=1}^{N+1}(D)_{jr}(u_{jr}^*+u_{jr}^n)+\frac{Pe}{2} \bigg( u_{j-1}^*+u_{j-1}^n-2 (u_j^*+u_j^n)+u_{j+1}^*+u_{j+1}^n \bigg) + \frac{Da}{2} (u^*_j+u^n_j), \label{appendB:uN+1}
\end{alignat}
\label{appendB:u}
\eseq
where $D$ is the matrix of the OUCS3 scheme given by equations~\eqref{eqn:CompStencil}, \eqref{appendB:interior}, \eqref{appendB:boundary}. The numerical amplification factor, $G_{\text{num}} = \nicefrac{u^{n+1}_j}{u^n_j}$, is given by,
\beq
G_{\text{num}}=1-\bigg[\sum^{N+1}_{r=1} \frac{N_c}{2} (D)_{jr} e^{ikh(r-j)}-Pe\{\cos(kh)-1\}-\frac{D_a}{2}\bigg] (1+G^*_{\text{num}}), \label{appendB:G}
\eeq
where
\beq
G^*_{\text{num}}=1-\sum^{N+1}_{r=1} N_c (D)_{jr} e^{ikh(r-j)}+2Pe\{\cos(kh)-1\}+D_a. \label{appendB:G*}
\eeq

\begin{figure}[htbp]
\centering
\begin{subfigure}{0.48\textwidth}
\includegraphics[width=0.96\linewidth, height=0.8\linewidth]{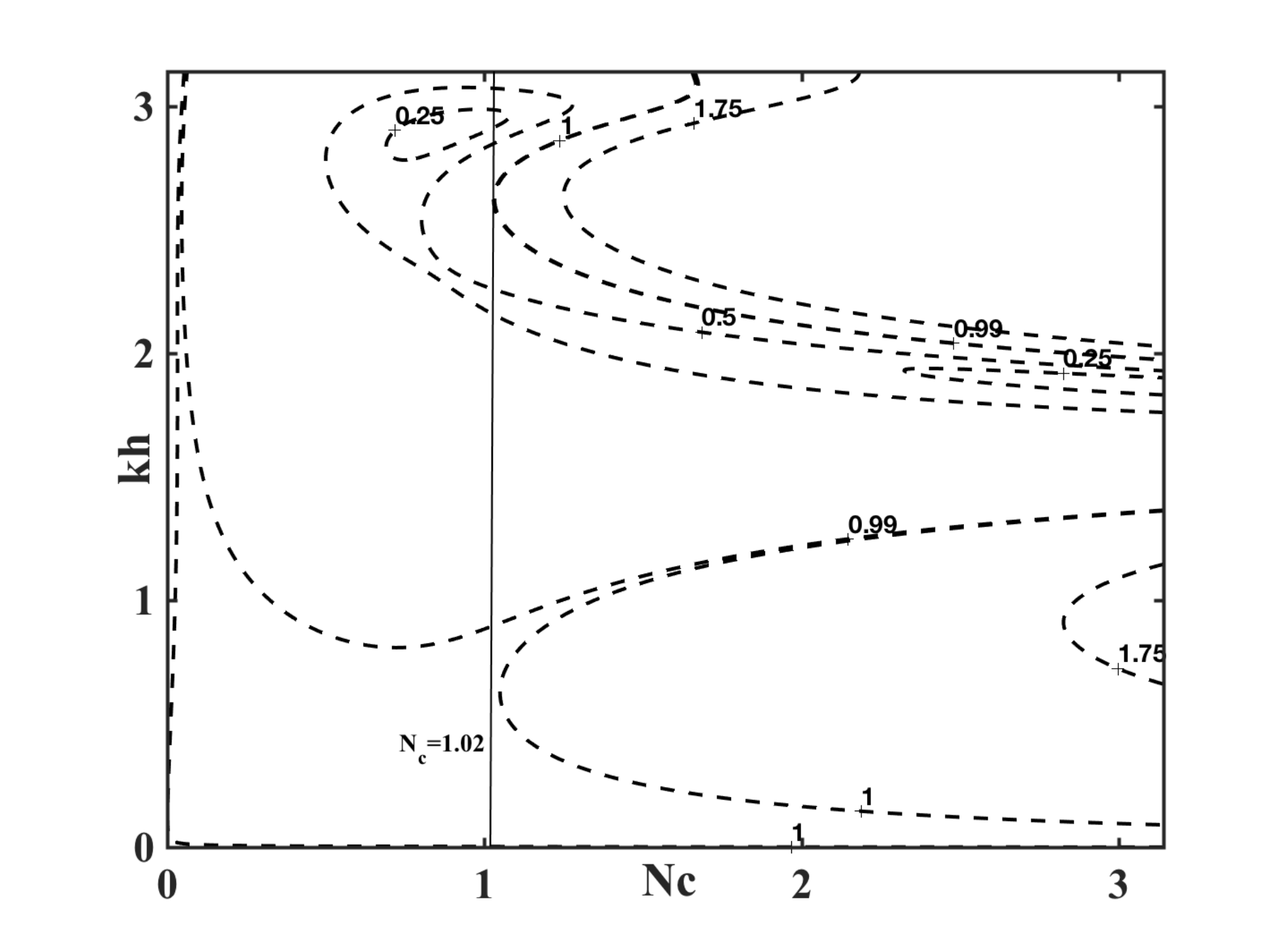}
\caption{} \label{fig:Fig1a}
\end{subfigure}
\begin{subfigure}{0.48\textwidth}
 \includegraphics[width=0.96\linewidth, height=0.8\linewidth]{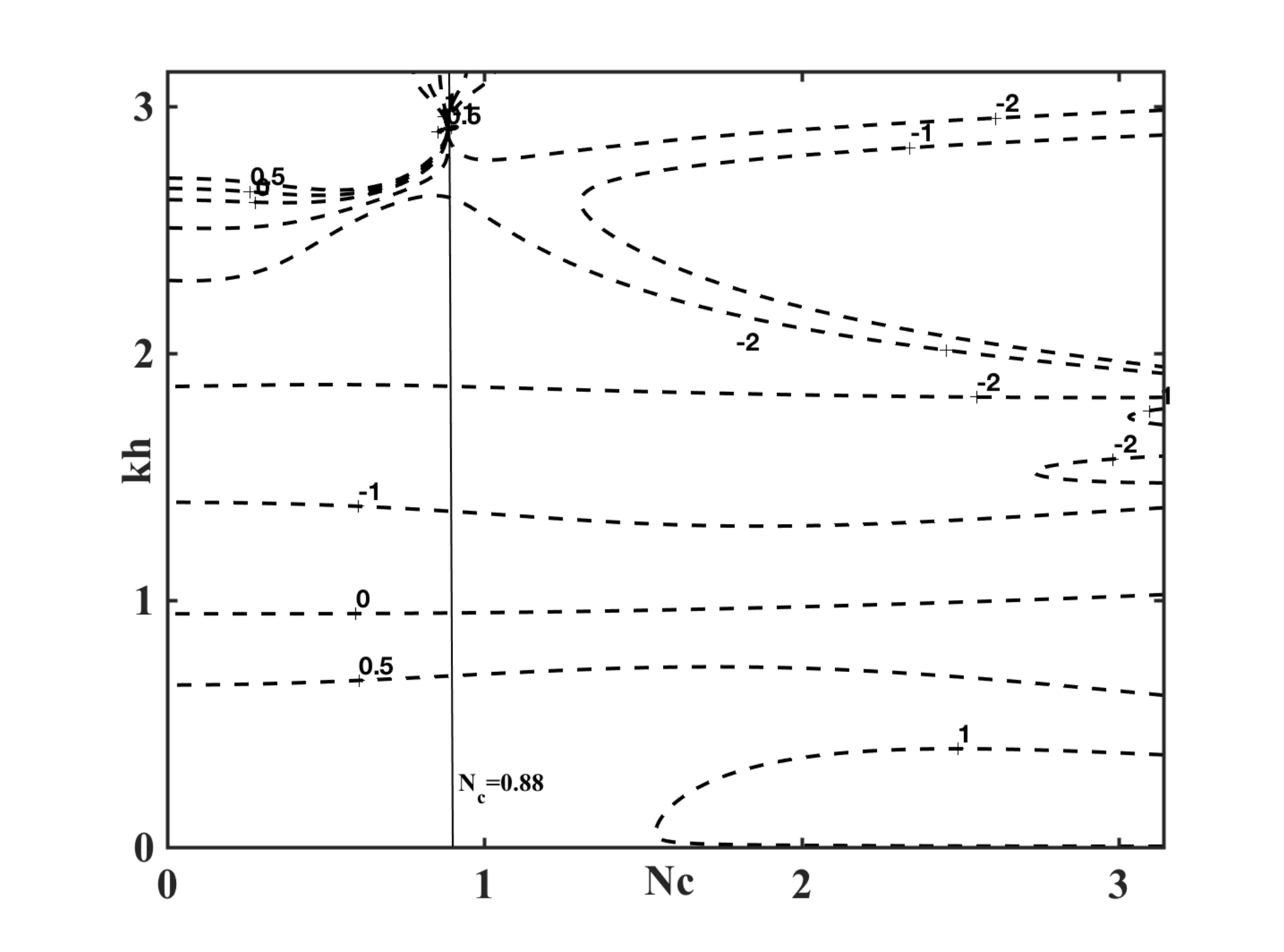}
 \caption{} \label{fig:Fig1b}
\end{subfigure}
\begin{subfigure}{0.48\textwidth}
 \includegraphics[width=.96\linewidth, height=0.8\linewidth]{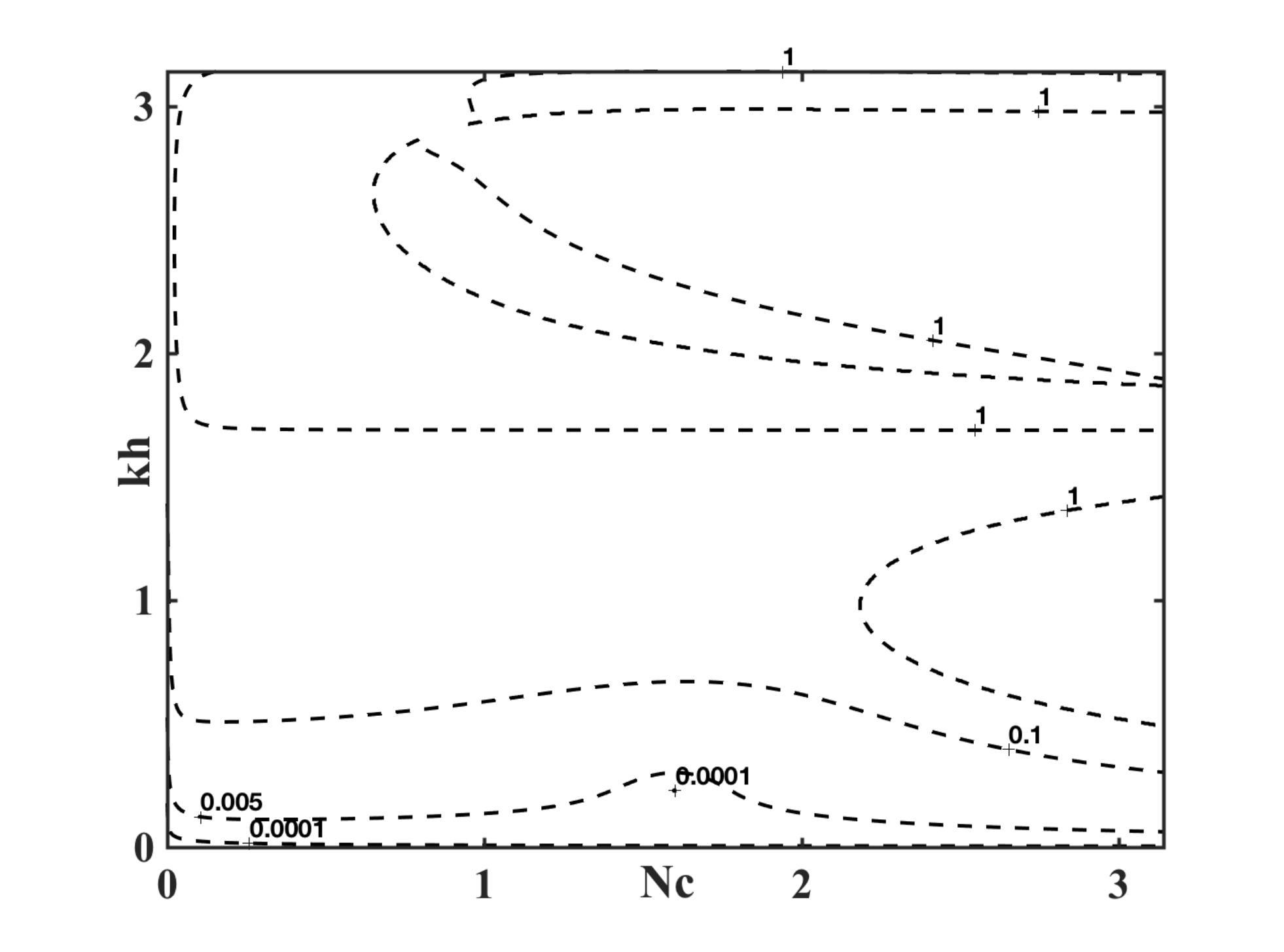}
 \caption{} \label{fig:Fig1c}
\end{subfigure}
\caption{Contour plots showing the (a) absolute amplification factor ratio, $\nicefrac{G_{\text{num}}}{G_{\text{exact}}}$, (b) scaled group velocity, $\nicefrac{V_{g,\text{num}}}{V_{g,\text{exact}}}$, and (c) absolute phase error, $|1 - \nicefrac{c_{\text{num}}}{c_{\text{exact}}}|$, using RK$_2$/Heun's time marching and OUCS3/CD$_2$ spatial discretization, (\S \ref{subsec:scheme1}, equations~(\ref{appendB:interior}-\ref{appendB:G})).}
\label{fig:Fig1}
\end{figure}

In figure~\ref{fig:Fig1a}, the amplification ratio, $G$, are plotted as contours in the indicated range of $N_c$ and $kh$. One notes that in the limit of vanishingly small values of $N_c$ and $kh$, the stability of this numerical scheme is reaction-rate dependent (i.~e., in the limit, $kh \rightarrow 0$ and $N_c \rightarrow 0$, we have $\lim\,\,G \rightarrow (1+0.5 Da)(1+Da)e^{-Da}$, refer equations~(\ref{appendB:G}, \ref{appendB:G*}, \ref{appendA:Gexact}). However, for intermediate values of the CFL number, $0 < N_c \le 1.02$, this numerical method is stable (or $G \le 1$). In figure~\ref{fig:Fig1b}, the scaled group velocity contours show significant dispersion effects, for almost all the selected values of $kh$ and $N_c$ (i.~e., $V_g \ne 1$), that would invalidate the long time integration results even when the neutral stability (i.~e., $G=1$) is ensured for those wavenumbers and CFL numbers. In fact the numerical solution travels in the wrong direction (or $V_g < 0$) within the range $N_c < 0.88,\,\,1.0 < kh \le 2.67$ or within $N_c \ge 0.88,\,\,1.0 \le kh \le \pi$ leading to numerical instabiltites in the form of $q-$waves (detailed analysis in \S \ref{sec:ED})~\cite{Vichnevetsky1982}. Although the error in the phase speed disappear for very small values of $kh$ (as highlighted in figure~\ref{fig:Fig1c}), we conclude through figures~\ref{fig:Fig1b}, \ref{fig:Fig1c} that the combined effects of the dispersion errors cannot be simply eliminated or reduced by grid refinement, as suggested in~\cite{Gustafsson1972}.

%%%%%%%%%%%%%%%%%%%%%%%%%%%%%%%%%%%%%%%%%%%%%%%%
\subsection{Implicit-OUCS3-Lele} \label{subsec:scheme2}
In the second case, equation~\eqref{eqn:ADR} is numerically time-integrated using the implicit mid-point rule~\cite{Stern2009}. The advection and the diffusion terms in equation~\eqref{eqn:ADR} are spatially discretized using the OUCS3 scheme, equations~(\ref{appendB:interior}, \ref{appendB:boundary}) and the Lele's scheme, respectively. The interior, boundary and the near boundary stencils of the Lele's scheme for the second derivative~\cite{Lele1992} are given by,
\bseq
\begin{alignat}{4}
&j = 1: \quad u''_1 = \frac{1}{h^2}(u_1 - 2 u_2 + u_3), \label{appendC:upp1} \\
&j = 2: \quad u''_1 + 10 u''_2 + u''_3 = \frac{12}{h^2}(u_3 - 2u_2 + u_1), \label{appendC:upp2} \\
&j \in [3, N-1]: \quad \alpha u''_{j-1} + u''_j + \alpha u''_{j+1} = \frac{b}{4 h^2}(u_{j-2} - 2 u_j + u_{j+2}) + \frac{a}{h^2}(u_{j-1} - 2 u_j + u_{j+1}), \label{appendC:interior} \\
&j = N: \quad u''_{N+1} + 10 u''_N + u''_{N-1} = \frac{12}{h^2}(u_{N-1} - 2 u_N + u_{N+1}),\label{appendC:uppN} \\
&j = N+1: \quad u''_{N+1} + 11 u''_N = \frac{1}{h^2} (13 u_{N+1} - 27 u_N + 15 u_{N-1} - u_{N-2}). \label{appendC:uppN+1}
\end{alignat}
\label{appendC:Lele}
\eseq
Lele's compact scheme~\cite{Lele1992} has superior resolution (e.~g.~when compared with the repeated operation of evaluating the first derivative twice using OUCS3) and has lower dispersion error for all the nodes~\cite{TKS2013}. The solution of equation~\eqref{eqn:ADR} at the $j$th node, $u^{n+1}_j$, is expressed as a function of the known at the previous time-step, $u^n_j$, as follows,
\beq
\bigg[\bigg( 1-\frac{Da}{2} \bigg)I +\frac{N_c}{2}(D_1)_{jr}-\frac{Pe}{2}(D_2)_{jr} \bigg]u^{n+1}_{jr}= 
\bigg[\bigg( 1+\frac{Da}{2} \bigg)I -\frac{N_c}{2}(D_1)_{jr}+\frac{Pe}{2}(D_2)_{jr} \bigg]u^n_{jr},
\label{appendC:u}
\eeq
where $D_1$ and $D_2$ are the matrices of the OUCS3 scheme (equation~(\ref{appendB:interior}, \ref{appendB:boundary})) and the Lele's scheme (equation~\eqref{appendC:Lele}), respectively. The numerical amplification factor is specified as follows,
\beq
G_{\text{num}}=\frac{1+\frac{Da}{2}-\frac{1}{2}\sum^{N+1}_{r=1}[Nc(D_1)_{jr}-Pe(D_2)_{jr}]e^{ikh(r-j)}}{1-\frac{Da}{2}+\frac{1}{2}\sum^{N+1}_{r=1}[Nc(D_1)_{jr}-Pe(D_2)_{jr}]e^{ikh(r-j)}}.
\label{appendC:G}
\eeq
% equation~\eqref{appendC:Lele}, respectively (refer \S \ref{sec:appendC} for details). 

\begin{figure}[htbp]
\centering
\begin{subfigure}{0.48\textwidth}
\includegraphics[width=0.96\linewidth, height=0.8\linewidth]{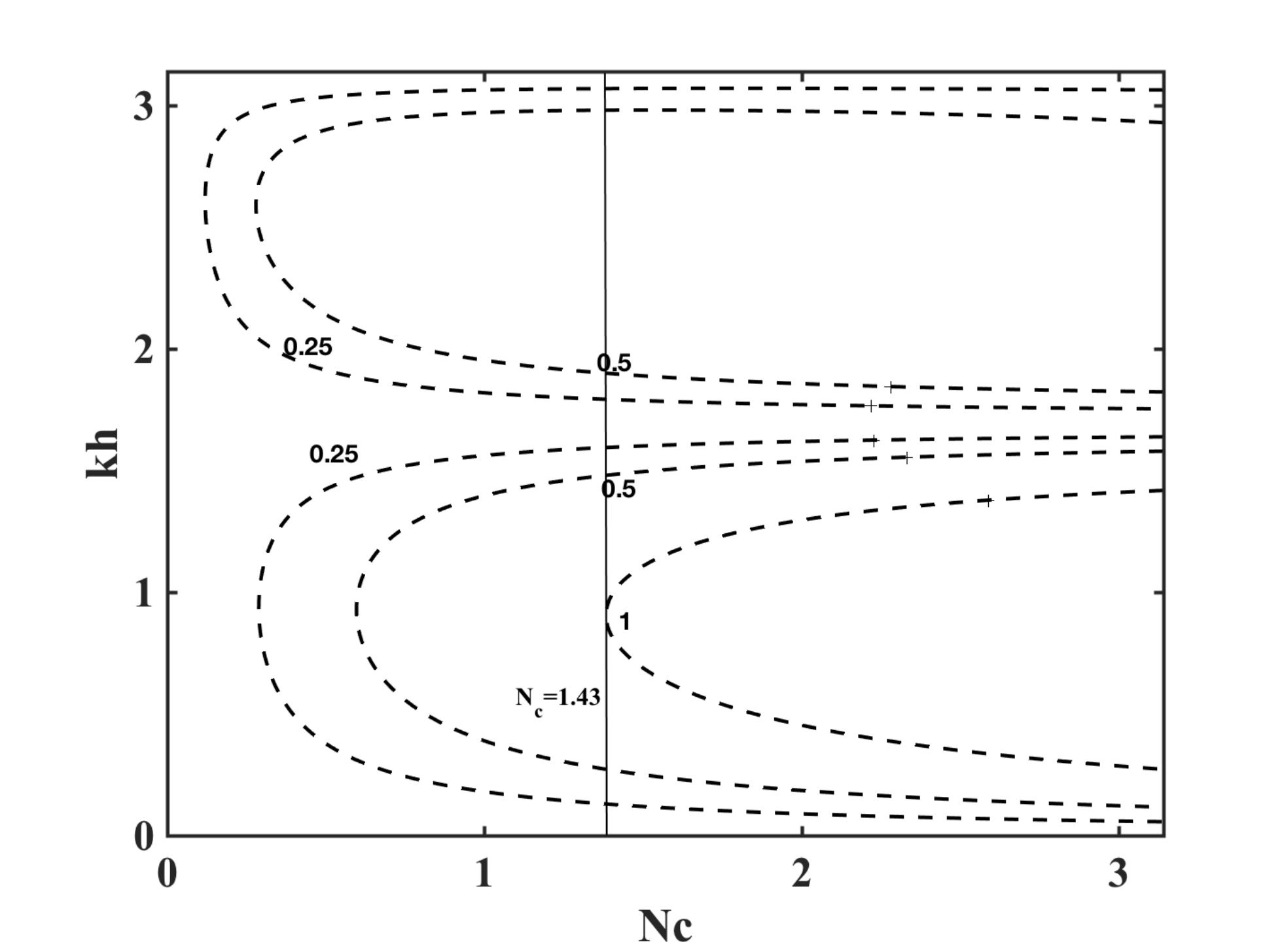}
\caption{} \label{fig:Fig2a}
\end{subfigure}
\begin{subfigure}{0.48\textwidth}
 \includegraphics[width=0.96\linewidth, height=0.8\linewidth]{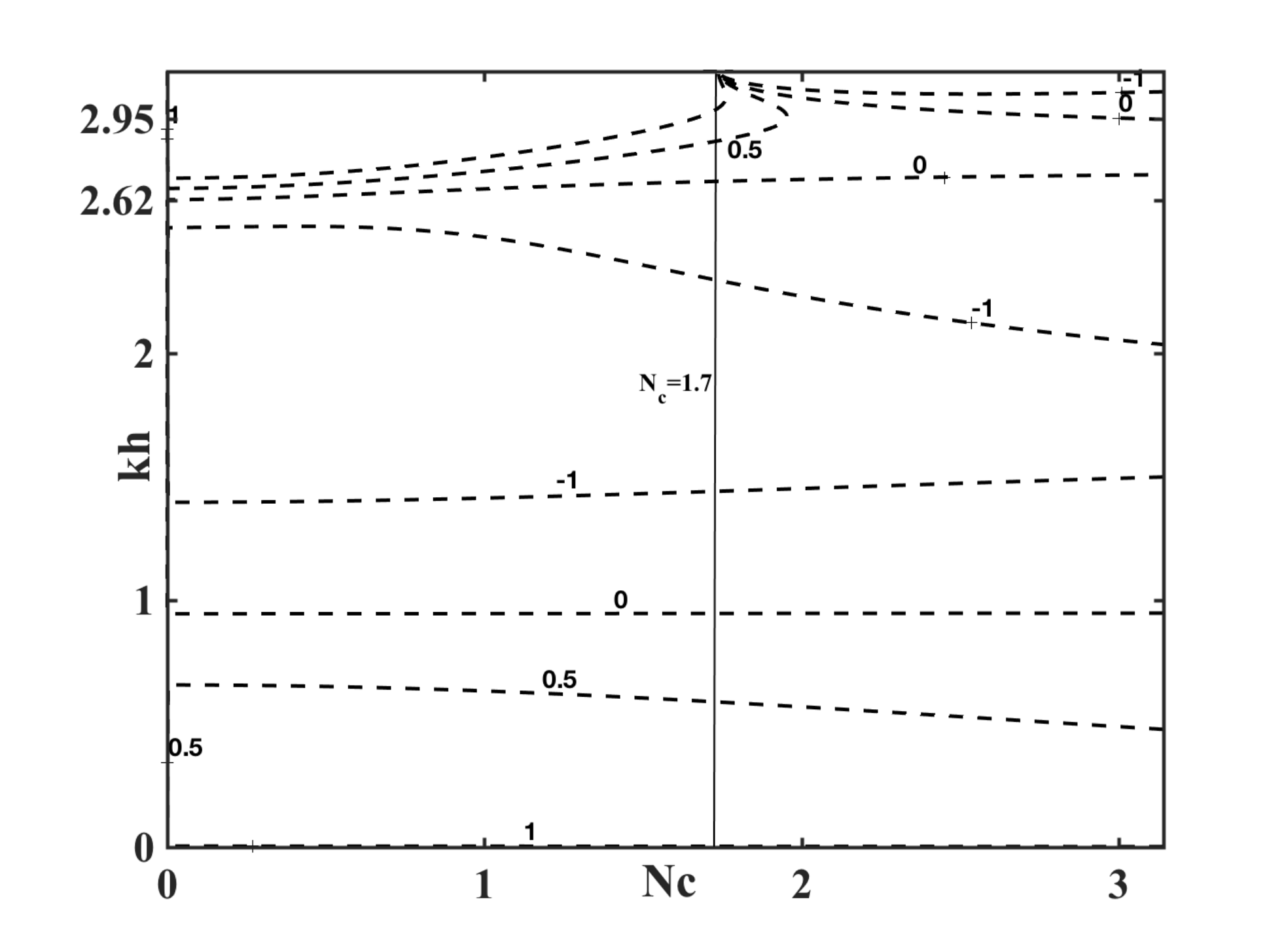}
 \caption{} \label{fig:Fig2b}
\end{subfigure}
\begin{subfigure}{0.48\textwidth}
 \includegraphics[width=.96\linewidth, height=0.8\linewidth]{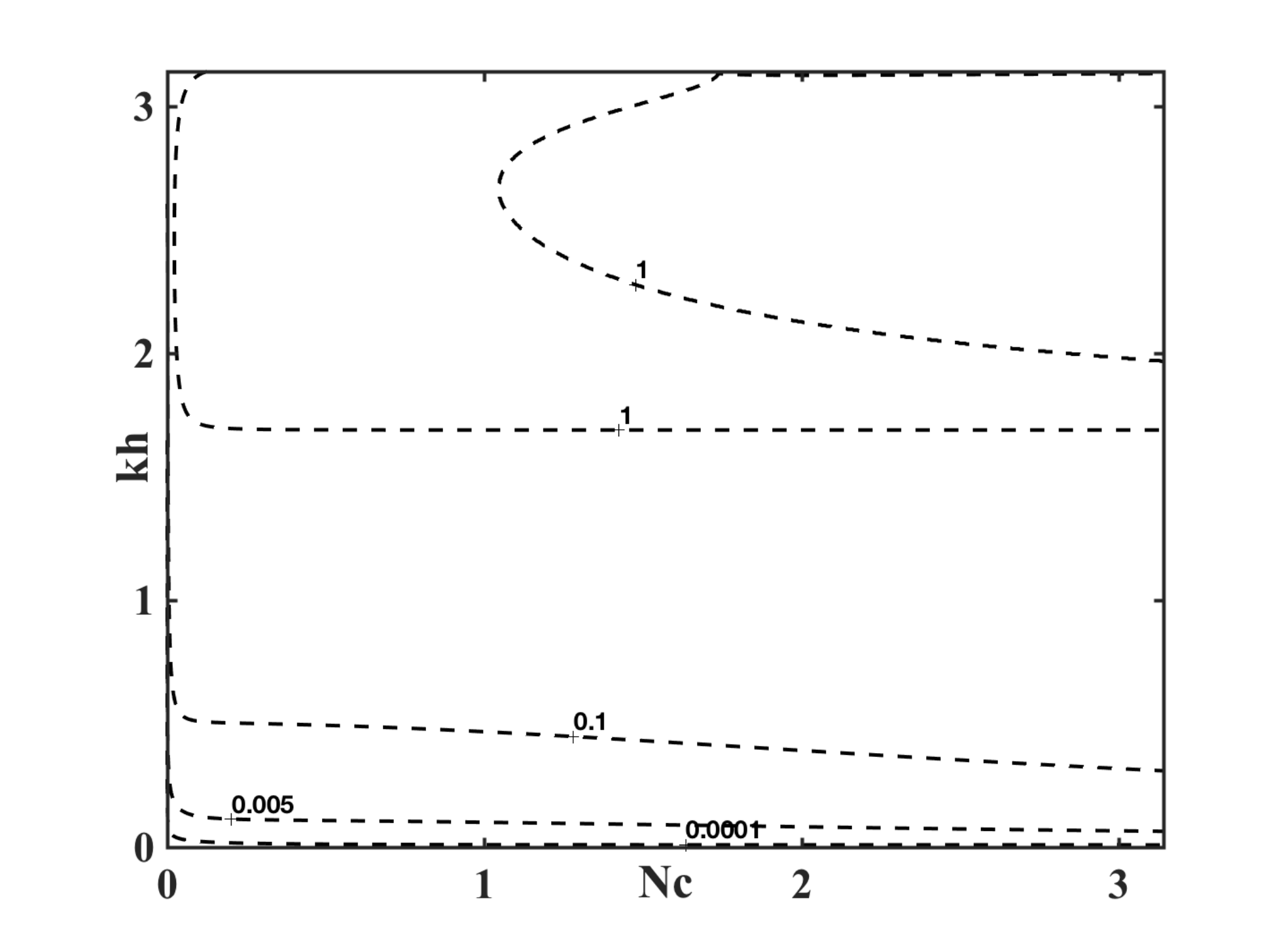}
 \caption{} \label{fig:Fig2c}
\end{subfigure}
\caption{Contour plots highlighting the (a) absolute amplification factor ratio, $\nicefrac{G_{\text{num}}}{G_{\text{exact}}}$, (b) scaled group velocity, $\nicefrac{V_{g,\text{num}}}{V_{g,\text{exact}}}$, and (c) absolute phase error, $|1 - \nicefrac{c_{\text{num}}}{c_{\text{exact}}}|$, using implicit mid-point time marching and OUCS3/Lele spatial discretization, (\S \ref{subsec:scheme2}, equations~(\ref{appendB:interior}-\ref{appendB:G}, \ref{appendC:Lele})).}\label{fig:Fig2}
\end{figure}

The contours of the amplification ratio are outlined in figure~\ref{fig:Fig2a} in the selected range of $kh$ and $N_c$. This discretization scheme has a comparatively larger stability region than the explicit scheme discussed in \S \ref{subsec:scheme1} (i.~e., $G \le 1$, for all $kh$ and $N_c \le 1.43$), a feature attributed due to the favorable temporal stability property of implicit methods for stiff equations~\cite{Stern2009}. The scaled group velocity contours in figure~\ref{fig:Fig2b} indicate that the region of spurious propagation (via $q-$waves) is limited within the range $1.0 < kh < 2.62$ (for all $N_c$) or within $kh > 2.95$ (and $N_c > 1.70$). Further, the dispersion errors via the scaled group velocity is inconsequential in the limit of small wavenumbers (i.~e., V$_g \rightarrow 1$ in the limiting case of $kh \rightarrow 0$). Figure~\ref{fig:Fig2c} suggests that the error in the phase speed vanishes in limit of negligibly small wavenumbers (i.~e., $c^{\text{err}}_{\text{phase}} \rightarrow 0$ when $kh \rightarrow 0$), an observation parallel to the numerical method discussed in \S \ref{subsec:scheme1} (refer figure~\ref{fig:Fig1c}). Overall, the stability and the propagation characteristics indicate that this discretization scheme is marginally better than the first case in determining the numerical solution of the 1-D linear advection-diffusion-reaction equation~\eqref{eqn:ADR}.
%, albeit the unwanted dispersion due to the propagation error is present in almost the entire region of $kh$-$N_c$ parameter space.

%%%%%%%%%%%%%%%%%%%%%%%%%%%%%%%%%%%%%%%%%%%%%%%%
\subsection{IMEX-OUCS3-Lele} \label{subsec:scheme3}
In the third case, the temporal discretization of equation~\eqref{eqn:ADR} is accomplished via a combination of the explicit two-stage Runge Kutta/Heun's method and the implicit mid-point rule~\cite{Stern2009}; the (nonstiff) advection term is treated explicitly while the (stiff) diffusion and reaction terms are handled implicitly. Spatial discretization is achieved through the OUCS3 (Lele) for first (second) derivatives, i.~e.~equation~(\ref{appendB:interior}, \ref{appendB:boundary}) (equation~(\ref{appendC:Lele})). Using this combination of spatiotemporal discretization, the unknown solution of equation~\eqref{eqn:ADR} at the $j$th node, $u^*_j, u^{n+1}_j$, at each stage of the temporal discretization as a function of the known variable at the previous time-step, $u^n_j$, is described as follows,
\bseq
\begin{alignat}{4}
&\bigg[ \bigg( 1-\frac{Da}{2}\bigg) I -\frac{Pe}{2} \sum^{N+1}_{r=1} (D_2)_{jr} \bigg] u_{jr}^*= \bigg[ \bigg(1+\frac{Da}{2}\bigg) I + \frac{Pe}{2} \sum^{N+1}_{r=1} (D_2)_{jr} - N_c \sum^{N+1}_{r=1} (D_1)_{jr} \bigg] u_{jr}^n, \label{appendD:u*}
\\
&u_j^{n+1}=u^n_j - \frac{1}{2} \bigg( N_c \sum^{N+1}_{r=1} (D_1)_{jr} - Pe \sum^{N+1}_{r=1} (D_2)_{jr} - (Da) I \bigg)(u^{n}_{jr} + u^*_{jr}), \label{appendD:uN+1}
\end{alignat}
\label{appendD:u}
\eseq
where $D_1$ and $D_2$ are the matrices of the OUCS3 scheme and the Lele's scheme, respectively. The numerical amplification factor is prescribed by,
\beq
G_{\text{num}}=1-\bigg(\sum^{N+1}_{r=1} \bigg[\frac{N_c}{2}(D_1)_{jr}-\frac{Pe}{2}(D_2)_{jr} \bigg]\,\,e^{ikh(r-j)}-\frac{Da}{2}\bigg)(1+G^*_{\text{num}}), \label{appendD:G}
\eeq
where
\beq
G^*_{\text{num}}=1+\frac{Da-\sum^{N+1}_{r=1}[N_c (D_1)_{jr} - Pe (D_2)_{jr}]\, \, e^{ikh(r-j)}}{1-\frac{Da}{2}-\frac{Pe}{2}\sum^{N+1}_{r=1}(D_2)_{jr}\, \, e^{ikh(r-j)}}. \label{appendD:G*}
\eeq

\begin{figure}[htbp]
\centering
\begin{subfigure}{0.48\textwidth}
\includegraphics[width=0.96\linewidth, height=0.8\linewidth]{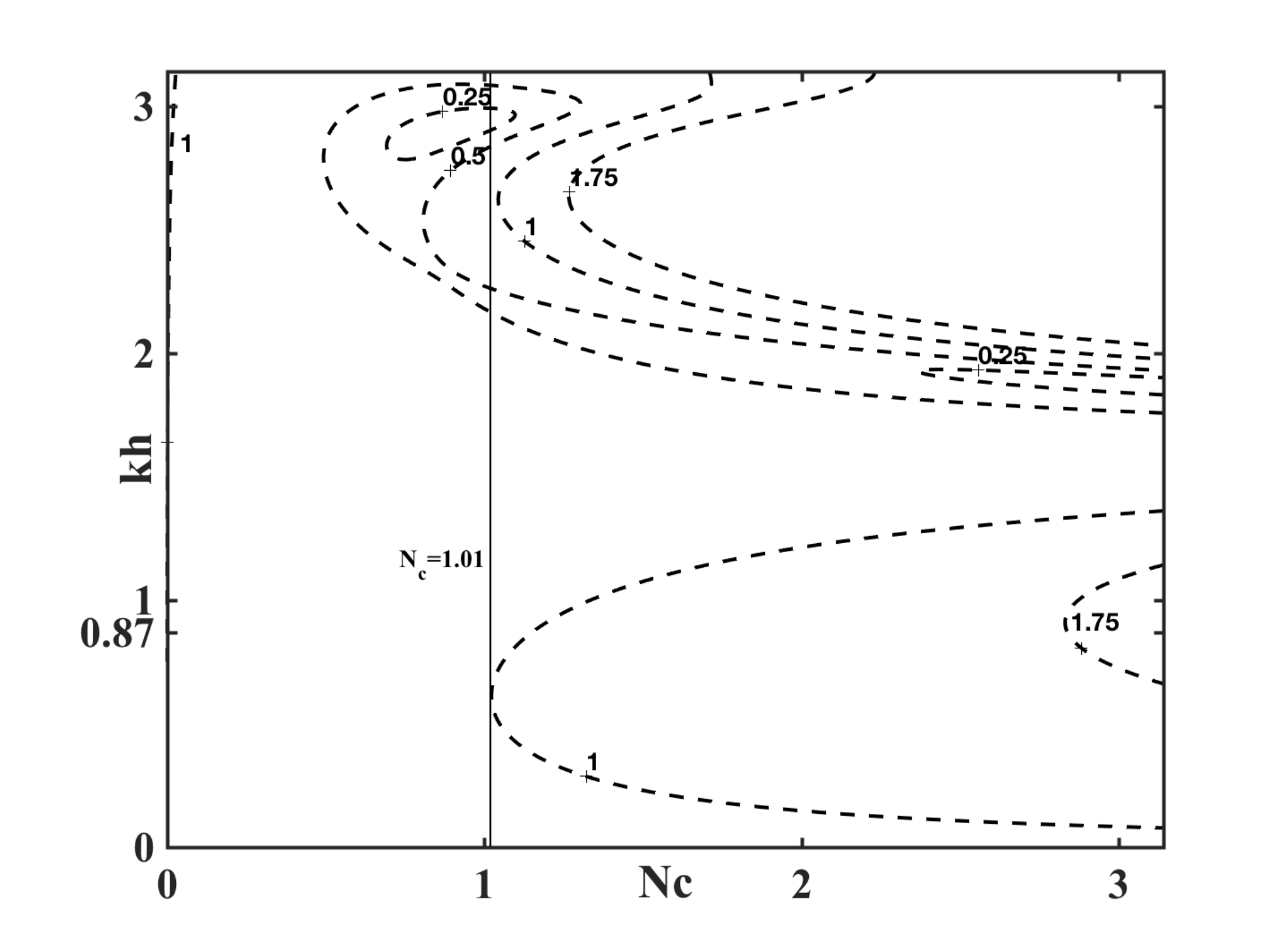}
\caption{} \label{fig:Fig3a}
\end{subfigure}
\begin{subfigure}{0.48\textwidth}
 \includegraphics[width=0.96\linewidth, height=0.8\linewidth]{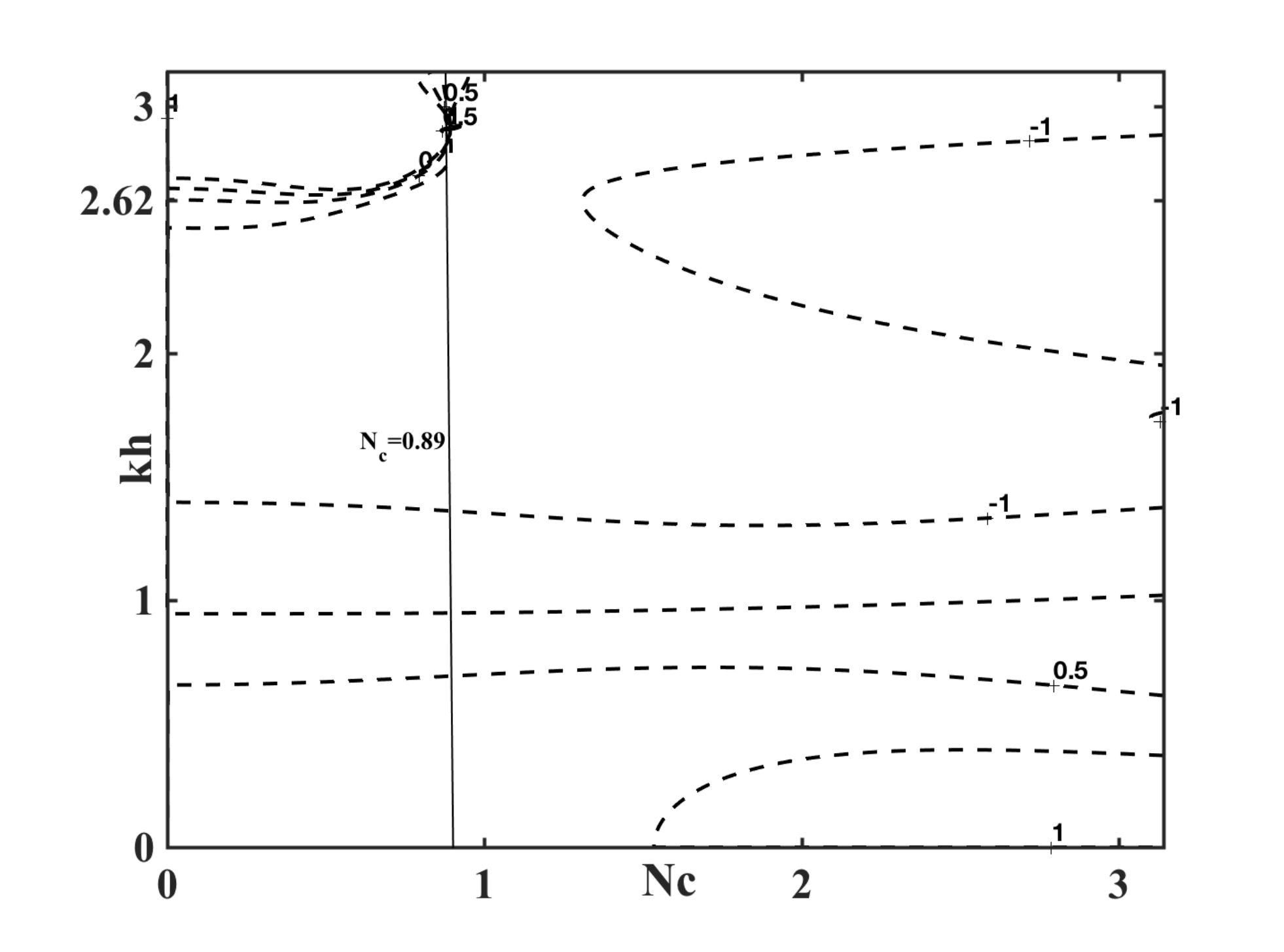}
 \caption{} \label{fig:Fig3b}
\end{subfigure}
\begin{subfigure}{0.48\textwidth}
 \includegraphics[width=.96\linewidth, height=0.8\linewidth]{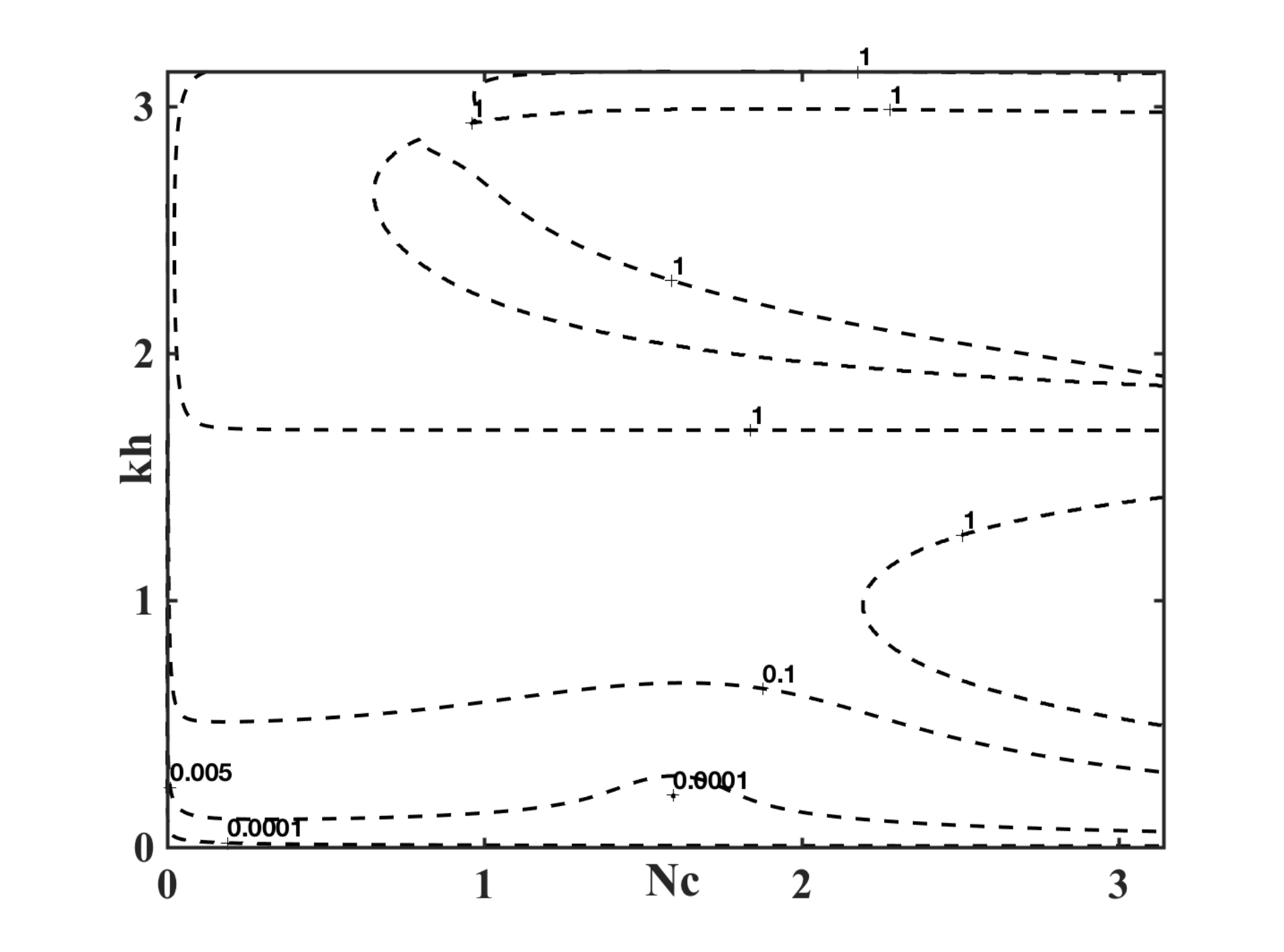}
 \caption{} \label{fig:Fig3c}
\end{subfigure}
\caption{Contour plots presenting the (a) absolute amplification factor ratio, $\nicefrac{G_{\text{num}}}{G_{\text{exact}}}$, (b) scaled group velocity, $\nicefrac{V_{g,\text{num}}}{V_{g,\text{exact}}}$, and (c) absolute phase error, $|1 - \nicefrac{c_{\text{num}}}{c_{\text{exact}}}|$, using IMEX (mid-point/RK$_2$) time marching and OUCS3/Lele spatial discretization, (\S \ref{subsec:scheme3}, equations~(\ref{appendB:interior}-\ref{appendB:G}, \ref{appendC:Lele})).}\label{fig:Fig3}
\end{figure}

The amplification ratio contours in this case indicate a region of stability, $N_c \le 1.01$ and $kh \le 0.87$ (figure~\ref{fig:Fig3a}), while the spurious propagation features (i.~e., $V_g<0$, figure~\ref{fig:Fig3b}) is restricted inside the domain, $N_c < 0.89, 1.0 < kh \le 2.62$ or within the region, $N_c \ge 0.89, 1.0 \le kh \le \pi$. As expected, the phase speed error in this case is inconsequential for small wavenumbers ($c^{\text{err}}_{\text{phase}} \approx 0$ in the limit $kh \rightarrow 0$, figure~\ref{fig:Fig3c}). In summary, this numerical method is comparable to the second case (\S \ref{subsec:scheme2}) but with slightly reduced CFL range for temporal stability.
%marginally inferior stability property.
%

%%%%%%%%%%%%%%%%%%%%%%%%%%%%%%%%%%%%%%%%%%%%%%%%
\subsection{IMEX-NCCD} \label{subsec:scheme4}
In this numerical method, the IMEX temporal discretization (refer \S \ref{subsec:scheme3}) combined with the NCCD scheme for the spatial discretization of the advection and the diffusion terms is utilized to analyze equation~\eqref{eqn:ADR}~\cite{TKS2009P1,TKS2009P2}. The NCCD scheme is used to simultaneously evaluate the first and the second spatial derivatives ($u'_j, u''_j$) from the following discrete equations for the boundary and the interior nodes in terms of the known values, $u_j$,
\bseq
\begin{alignat}{4}
j = 1: \quad &u'_1 + 2u'_2 - h u''_2 = \frac{1}{h}(-3.5u_1 + 4u_2 - 0.5u_3), \label{appendE:up1} \\
&h u''_1 + 5h u''_2 - 6 u'_2 = \frac{1}{h}(9u_1 - 12u_2 + 3u_3), \label{appendE:upp1} \\
j \in [2, N]: \quad &\frac{7}{16} (u'_{j+1} + u'_{j-1}) + u'_j - \frac{h}{16} (u''_{j+1} - u''_{j-1}) = \frac{15}{16 h} (u_{j+1} - u_{j-1}), \label{appendE:stencil1} \\
&\frac{9}{8 h} (u'_{j+1} - u'_{j-1}) + u''_j - \frac{1}{8} (u''_{j+1} + u''_{j-1}) = \frac{3}{h^2} (u_{j+1} -2u_j + u_{j-1}), \label{appendE:stencil2} \\
j = N+1: \quad &u'_{N+1} + 2u'_N + h u''_N = -\frac{1}{h}(-3.5u_{N+1} + 4u_N - 0.5u_{N-1}), \label{appendE:upN+1} \\
&h u''_{N+1} + 5h u''_N + 6 u'_N = \frac{1}{h}(9u_{N+1} - 12u_N + 3u_{N-1}). \label{appendE:uppN+1}
\end{alignat}
\label{appendE:interior}
\eseq
%
%Equation~\eqref{appendE:stencil1} is used for $j = 3$ to $(N-1)$, while equation~\eqref{appendE:stencil2} is used for $j = 2$ to $j=N$. 
To avoid the instability and the attenuation near the inflow and the outflow, the first derivative of the unknowns near boundary ($u'_2, u'_N$) are eventually replaced with the locally explicit stencil~(\ref{appendB:up2}, \ref{appendB:upN}), while the second derivative of the unknowns near boundary ($u''_2, u''_N$) are replaced with~\cite{TKS2013},
\bseq
\begin{alignat}{4}
&u''_2 = \frac{1}{h^2}(u_1 - 2 u_2 + u_3), \label{appendE:upp1} \\
&u''_N = \frac{1}{h^2}(u_{N+1} - 2 u_N + u_{N-1}). \label{appendE:uppN+1} 
\end{alignat}
\label{appendE:boundary}
\eseq
Writing the NCCD interior and boundary stencils given by equations~(\ref{appendB:boundary}, \ref{appendE:interior}, \ref{appendE:boundary}) as a system of linear algebraic equations,
\bseq
\begin{alignat}{4}
&[A_1]\{ u' \} + [B_1]\{ u'' \} = [C_1]\{ u \}, \label{appendE:system1} \\
&[A_2]\{ u' \} + [B_2]\{ u'' \} = [C_2]\{ u \}, \label{appendE:system2}
\end{alignat}
\label{appendE:system}
\eseq
and on solving equations~(\ref{appendE:system1}, \ref{appendE:system2}) simultaneously we arrive at,
\bseq
\begin{alignat}{4}
&\{ u' \} = \frac{1}{h} [D_1] \{ u \}, \label{appendE:up} \\
&\{ u'' \} = \frac{1}{h^2} [D_2] \{ u \}, \label{appendE:upp}
\end{alignat}
\label{appendE:Ssystem}
\eseq
where,
\bseq
\begin{alignat}{4}
&[D_1] = ([A_1] - [B_1][B_2]^{-1}[A_2])^{-1}([C_1] - [B_1][B_2]^{-1}[C_2]), \label{appendE:D1} \\
&[D_2] = ([B_2] - [A_2][A_1]^{-1}[B_1])^{-1}([C_2] - [A_2][A_1]^{-1}[C_1]). \label{appendE:D2}
\end{alignat}
\label{appendE:NCCDmatrix}
\eseq
The numerical amplification factor is evaluated using equations~(\ref{appendD:G}, \ref{appendD:G*}) with $D_1, D_2$ replaced with the matrices of the NCCD scheme given by equations~(\ref{appendE:D1},\ref{appendE:D2}), respectively.

\begin{figure}[htbp]
\centering
\begin{subfigure}{0.48\textwidth}
\includegraphics[width=0.96\linewidth, height=0.8\linewidth]{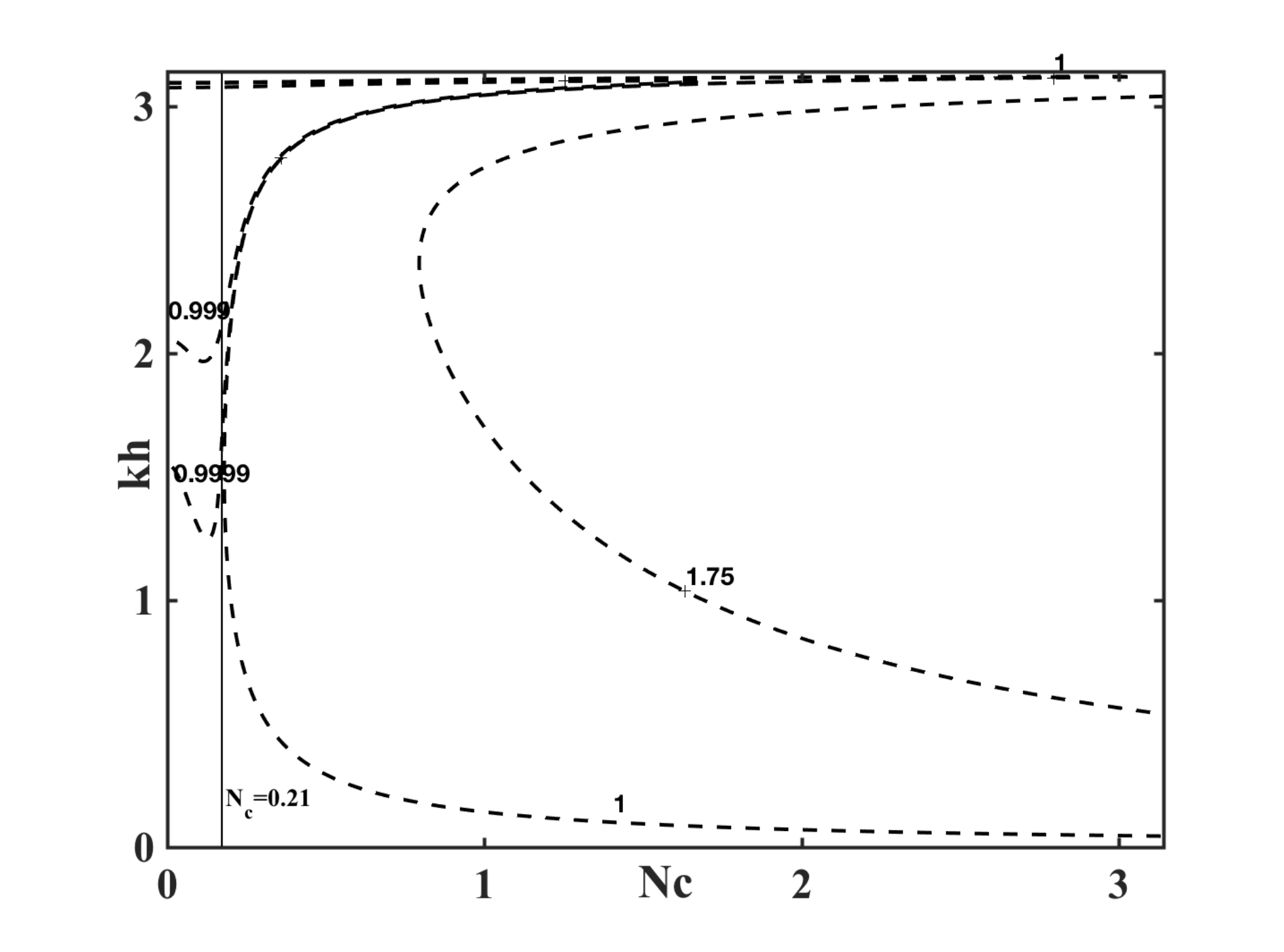}
\caption{} \label{fig:Fig4a}
\end{subfigure}
\begin{subfigure}{0.48\textwidth}
 \includegraphics[width=0.96\linewidth, height=0.8\linewidth]{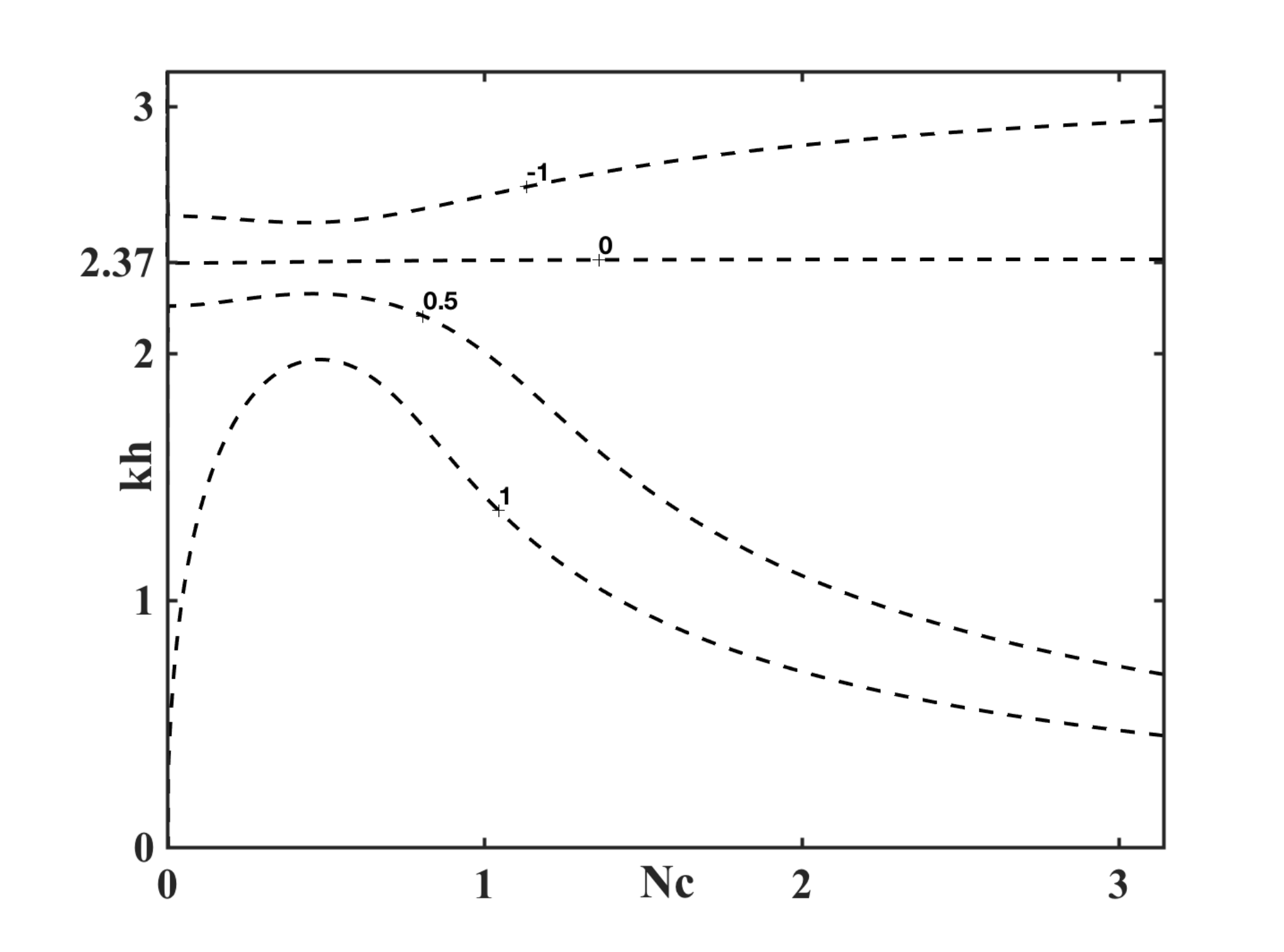}
 \caption{} \label{fig:Fig4b}
\end{subfigure}
\begin{subfigure}{0.48\textwidth}
 \includegraphics[width=.96\linewidth, height=0.8\linewidth]{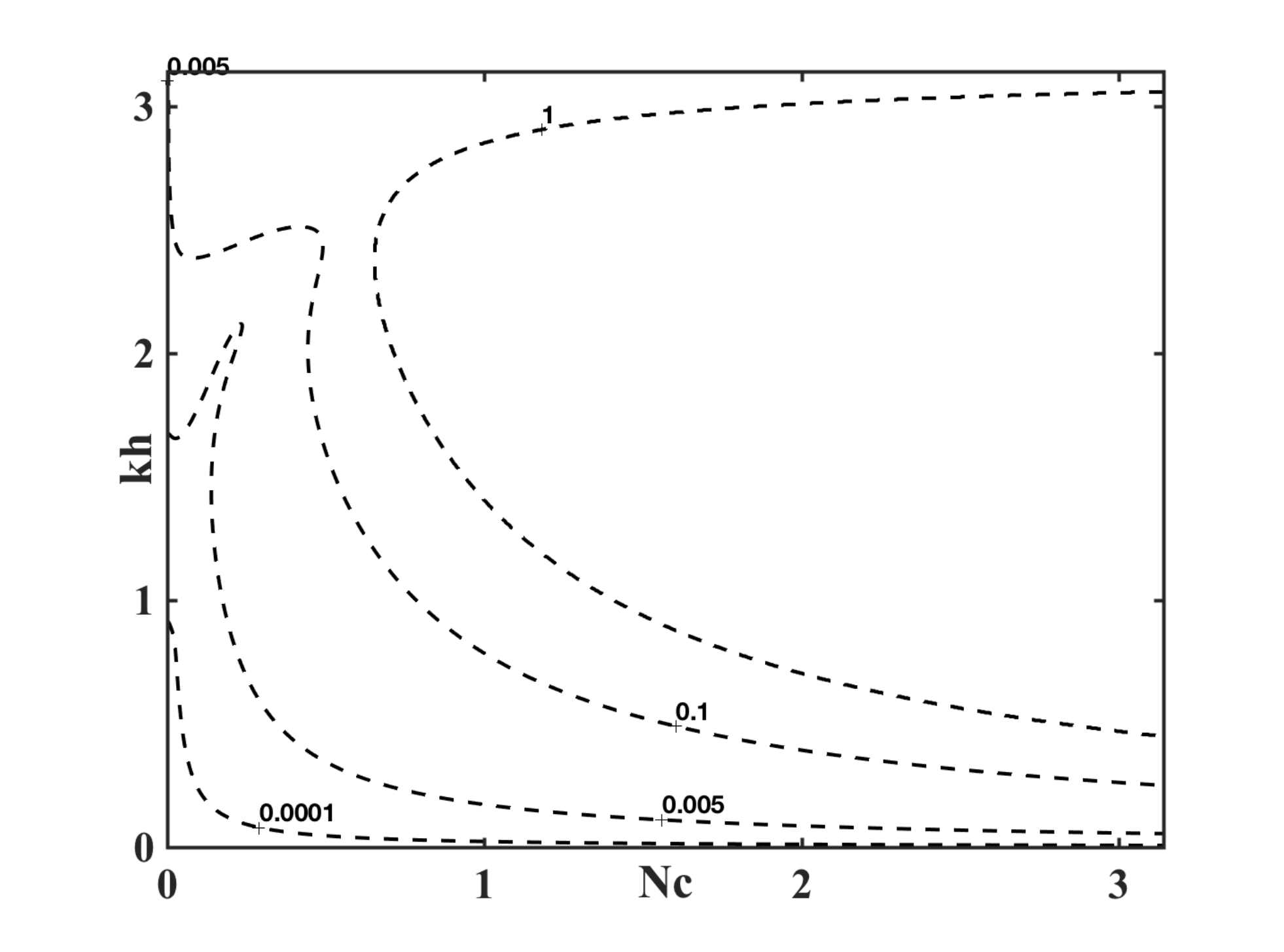}
 \caption{} \label{fig:Fig4c}
\end{subfigure}
\caption{Contour plots describing the (a) absolute amplification factor ratio, $\nicefrac{G_{\text{num}}}{G_{\text{exact}}}$, (b) scaled group velocity, $\nicefrac{V_{g,\text{num}}}{V_{g,\text{exact}}}$, and (c) absolute phase error, $|1 - \nicefrac{c_{\text{num}}}{c_{\text{exact}}}|$, using IMEX (mid-point/RK$_2$) time marching and NCCD spatial discretization, (\S \ref{subsec:scheme4}, equations~(\ref{appendB:boundary}, \ref{appendE:interior}-\ref{appendE:boundary})).}\label{fig:Fig4}
\end{figure}

The superior numerical resolution feature of this method is evident from the amplification ratio contours in figure~\ref{fig:Fig4a}, this scheme is neutrally stable (i.~e., $G \approx 1$) within the region, $N_c \le 0.21, 0 \le kh \le \pi$ - a property which is absolutely essential in Direct Numerical Simulation (DNS) of multiscale models~\cite{Canuto1987}. The spurious propagation characteristics (or the region of negative group velocity, figure~\ref{fig:Fig4b}) is limited to a relatively smaller range, $0 \le N_c \le \pi, kh > 2.37$, while the region where the dispersion effects due to phase speed error (figure~\ref{fig:Fig4c}) is diminished in comparison with the other three numerical methods (\S \ref{subsec:scheme1}-\ref{subsec:scheme3}). To recapitulate, the IMEX-NCCD scheme shows excellent DRP properties and it is undoubtedly the preferred numerical method to solve the 1D linear advection-diffusion-reaction equation~\eqref{eqn:ADR}, a conclusion which is further elucidated via the analysis through the error propagation equation, described next. 

%%%%%%%%%%%%%%%%%%%%%%%%%%%%%%%%%%%%%%%%%%%%%%%%
\section{Error dynamics} \label{sec:ED}
In this section, we identify the reasons for the dispersion error in the computational solution of equation~\eqref{eqn:ADR}, solved via the numerical methods discussed in \S \ref{sec:spectral}. The analyses in this section is based on the error propagation equation~\eqref{appendA:ErrorADR}. An extreme form of the dispersion error has been identified as the $q-$waves~\cite{Vichnevetsky1982} which are small amplitude, parasitic waves propagating in the direction opposite to the physical direction of propagation~\cite{TKS2012}, leading to numerical instabilities as well as unphysical bypass transition~\cite{TKS2010}. If $V_g$ is negative the numerical waves propagate upstream despite the physical requirement of downstream movement (since the physical advection speed, $c > 0$ in equation~\eqref{eqn:ADR}). In other words, $V_g < 0$ is the necessary condition for the generation of $q-$waves. Further elaboration on the role of $q-$waves in the solution of equation~\eqref{eqn:ADR} is considered through the propagation of a wave-packet which is given by the following initial condition,
\beq
u(x, t)|_{t=0} = e^{-\gamma(x-x_0)^2} \cos [k_0 (x - x_0)], \quad -L \le x \le L
\label{eqn:IC-ADR}
\eeq
where $\gamma$ controls the width of the packet and x$_0$, $k_0$ and $L$ is the center of the wave-packet, the central wavenumber and half domain length, respectively. The in silico analyses of equation~\eqref{eqn:ADR} together with the initial condition~\eqref{eqn:IC-ADR} and Dirichlet boundary condition ($u(\pm L, 0)$, equation~\eqref{eqn:IC-ADR}) is executed with the fixed parameters $c=0.1,  \nu=10^{-4}, \lambda=-1.0, x_0=0.0, k_0h=0.5, L=5.0$. The results reported in \S \ref{subsec:effect1}-\ref{subsec:effect3} are presented at (non-dimensional) simulation time $t=10.0$, using a time-step of numerical integration, $dt=0.01$ and (with the exception of the results described in \S \ref{subsec:effect2}) the domain is discretized uniformly with $N=1001$ points. Two trials cases are considered where the width of the initial wave-packet is (a) compact, $\gamma=50$, and (b) outstretched, $\gamma=10000$. The initial conditions (equation~\eqref{eqn:IC-ADR}) along with their Fourier spectra are displayed in figures~\ref{fig:Fig5a} and \ref{fig:Fig5b}, respectively.

%for the creation of $q-$waves in solving equation~\eqref{eqn:ADR} by comparing the dispersion error of the numerical methods discussed in \S \ref{sec:spectral}. $q-$waves are
%in terms of the scaled group velocity, $V_g$ (refer equation~\eqref{appendA:Vg}). 
%The $V_g$ contours (figures~\ref{fig:Fig1b}, \ref{fig:Fig2b}, \ref{fig:Fig3b}, \ref{fig:Fig4b}) indicate that the methods discussed in \S \ref{sec:spectral} produce dispersion error, except for infinitesimally small wavenumbers or small grid-spacing (i. e.~$V_g \rightarrow 1$ as $kh \rightarrow 0$). 

%%%%%%%%%%%%%%%%%%%%%%%%%%%%%%%%%%%%%%%%%%%%%%%%\input{JSC.synctex.gz}

\subsection{Role of initial conditions} \label{subsec:effect1}
The solution of equations~(\ref{eqn:ADR}, \ref{eqn:IC-ADR}) using the Explicit-OUCS3-CD$_2$ scheme (\S \ref{subsec:scheme1}) for $\gamma=50$ ($\gamma=10000$) is shown figure~\ref{fig:Fig5c} (figure~\ref{fig:Fig5d}). The reason behind the negligible presence of $q-$waves for the case of compact initial wave-packet (figure~\ref{fig:Fig5c}) versus a substaintial presence of the same for the outstretched initial wave-packet case (i. e., the boxed region, figure~\ref{fig:Fig5d}) can be explained via the respective Fourier transform of the initial condition (figure~\ref{fig:Fig5b}). From figure~\ref{fig:Fig1b}, we have noted the region of negative scaled group velocity for this numerical scheme to be $1.0 < kh \le 2.67$ (and $N_c < 0.88$). The Fourier spectra reveals that the amplitude of the outstretched initial wave-packet (dashed curve, figure~\ref{fig:Fig5b}) is significantly large within this range of wavenumbers, $kh$, and it is these components of initial conditions which are responsible for $q-$waves. The examples in figure~\ref{fig:Fig5c}, \ref{fig:Fig5d} are depicted for fixed value of the parameters, $k_0h=0.5, N_c=0.1$, and for these values the Explicit-OUCS3-CD$_2$ scheme predicts the scaled group velocity, $V_g=0.7078$. Thus, the numerical solution propagates at a speed slower than the exact solution while, for certain values of wavenumbers within the range $1.0 < kh \le 2.67$ (and $N_c < 0.88$), the $q-$waves travel at speed faster than the exact solution (refer figure~\ref{fig:Fig1b} for the contour values of $V_g$ within the wavenumber range $1.0 < kh \le 2.67$) leading to the oscillations as shown in figure~\ref{fig:Fig5d}.

%the presence of significantly large amplitude $q-$waves
%%%%%%%%%%%%%%%%%%%%%%%%%%%%%%%%%%%%%%%%%%%%%%%%
\subsection{Role of grid resolution} \label{subsec:effect2}
Spurious generation of flow structures in 1D advection equation, due to inadequate grid resolution leading to the creation and propagation of $q-$waves has been communicated in~\cite{TKS2012}. A similar analysis the 1D ADR equation is outlined in this section. Figure~\ref{fig:Fig5c}, \ref{fig:Fig5e} presents the solution of equations~(\ref{eqn:ADR}, \ref{eqn:IC-ADR}) using the Explicit scheme (\S \ref{subsec:scheme1}) for $\gamma=50$ with the domain discretized uniformly using $N=1001$ and $N=101$ points, respectively, such that the central wavenumber of the initial wave-packet $k_0 h = 0.5$ is fixed (equation~\eqref{eqn:IC-ADR}). Although the central wavenumber of the initial wave-packet is identical, the band-width of the waves resolved in both cases are different. A coarser grid (i.~e., $N=101$) is unable to capture the initial conditions accurately (solid curve, figure~\ref{fig:Fig5a}) resulting in a larger range of excited wavenumbers (as compared with the finer grid), in particular those wavenumbers that correspond to negative scaled group velocity thereby generating large amplitude $q-$waves. The in silico studies reported in~\cite{TKS2012} emphasize that merely introducing finer mesh is not adequate for accurate computing and one must be punctilious towards the basic DRP properties of the spatiotemporal discretization technique, as discussed in the next section.

%%%%%%%%%%%%%%%%%%%%%%%%%%%%%%%%%%%%%%%%%%%%%%%%
\subsection{Role of different numerical methods} \label{subsec:effect3}
\begin{figure}[htbp]
\centering
\begin{subfigure}{0.48\textwidth}
\includegraphics[width=0.96\linewidth, height=0.8\linewidth]{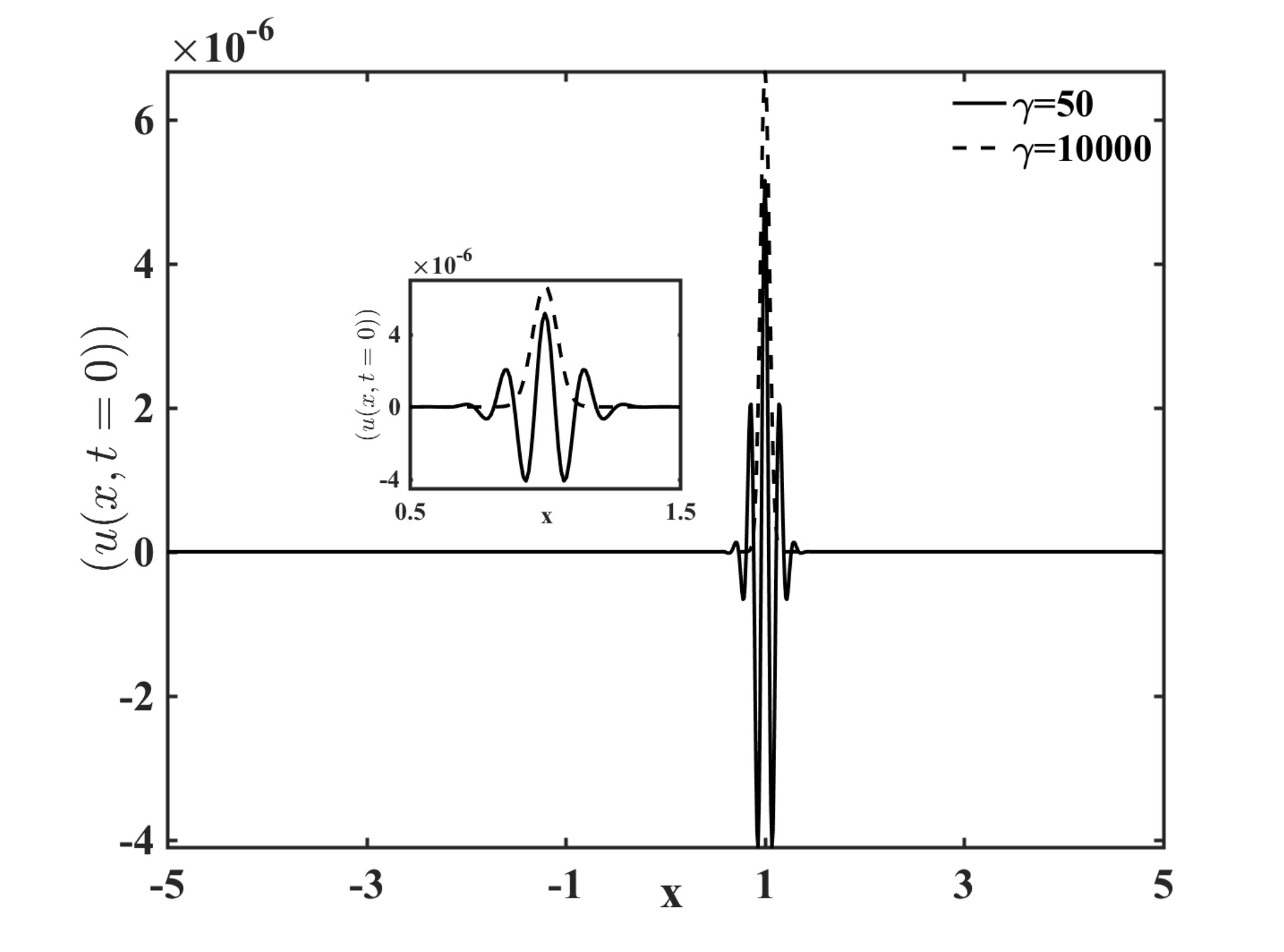}
\vskip -5pt
\caption{} \label{fig:Fig5a}
\end{subfigure}
\begin{subfigure}{0.48\textwidth}
 \includegraphics[width=0.96\linewidth, height=0.8\linewidth]{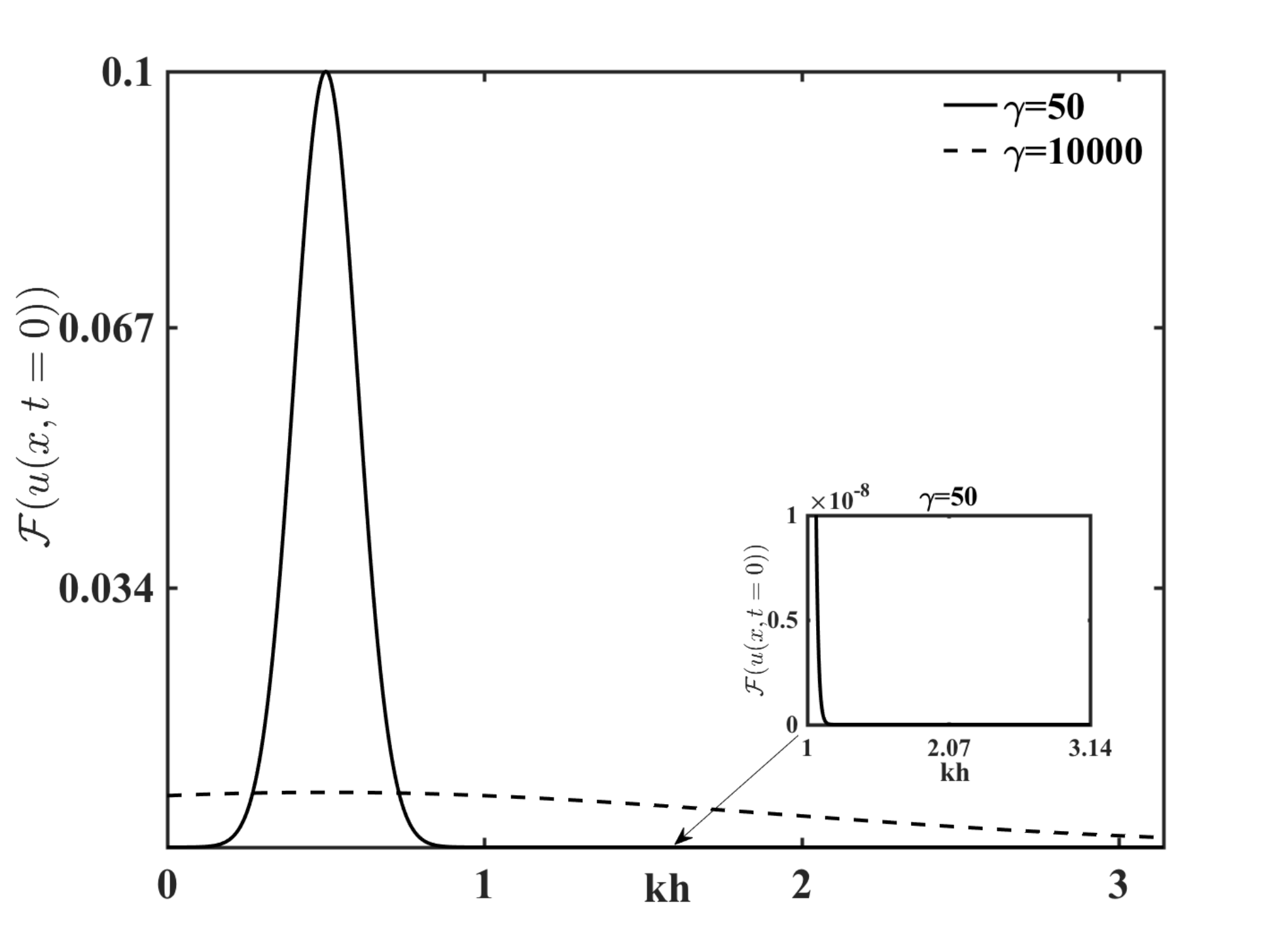}
 \vskip -5pt
 \caption{} \label{fig:Fig5b}
\end{subfigure}
\vskip -3.5pt
\begin{subfigure}{0.48\textwidth}
\includegraphics[width=0.96\linewidth, height=0.8\linewidth]{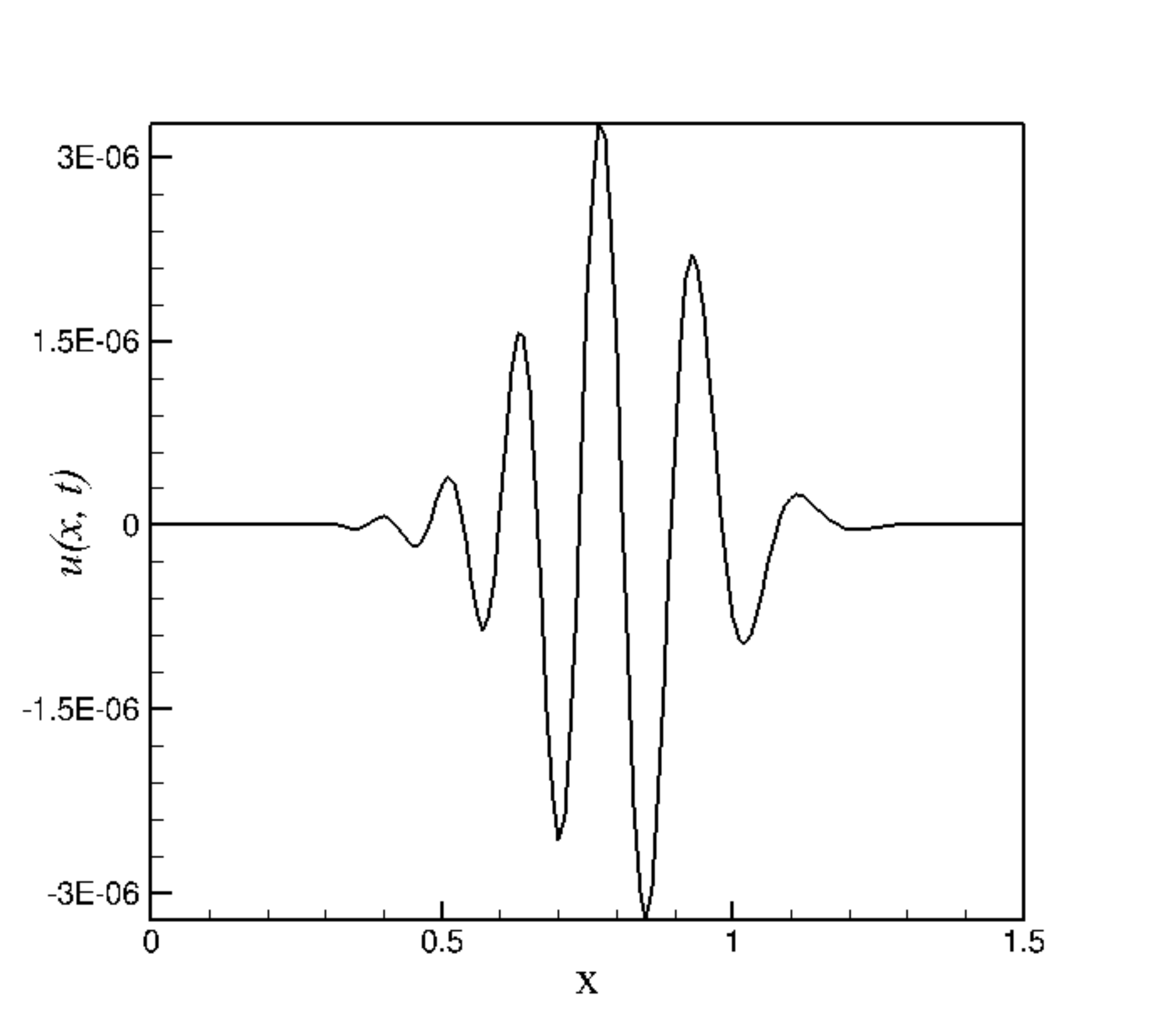}
\vskip -5pt
\caption{} \label{fig:Fig5c}
\end{subfigure}
\begin{subfigure}{0.48\textwidth}
 \includegraphics[width=0.96\linewidth, height=0.8\linewidth]{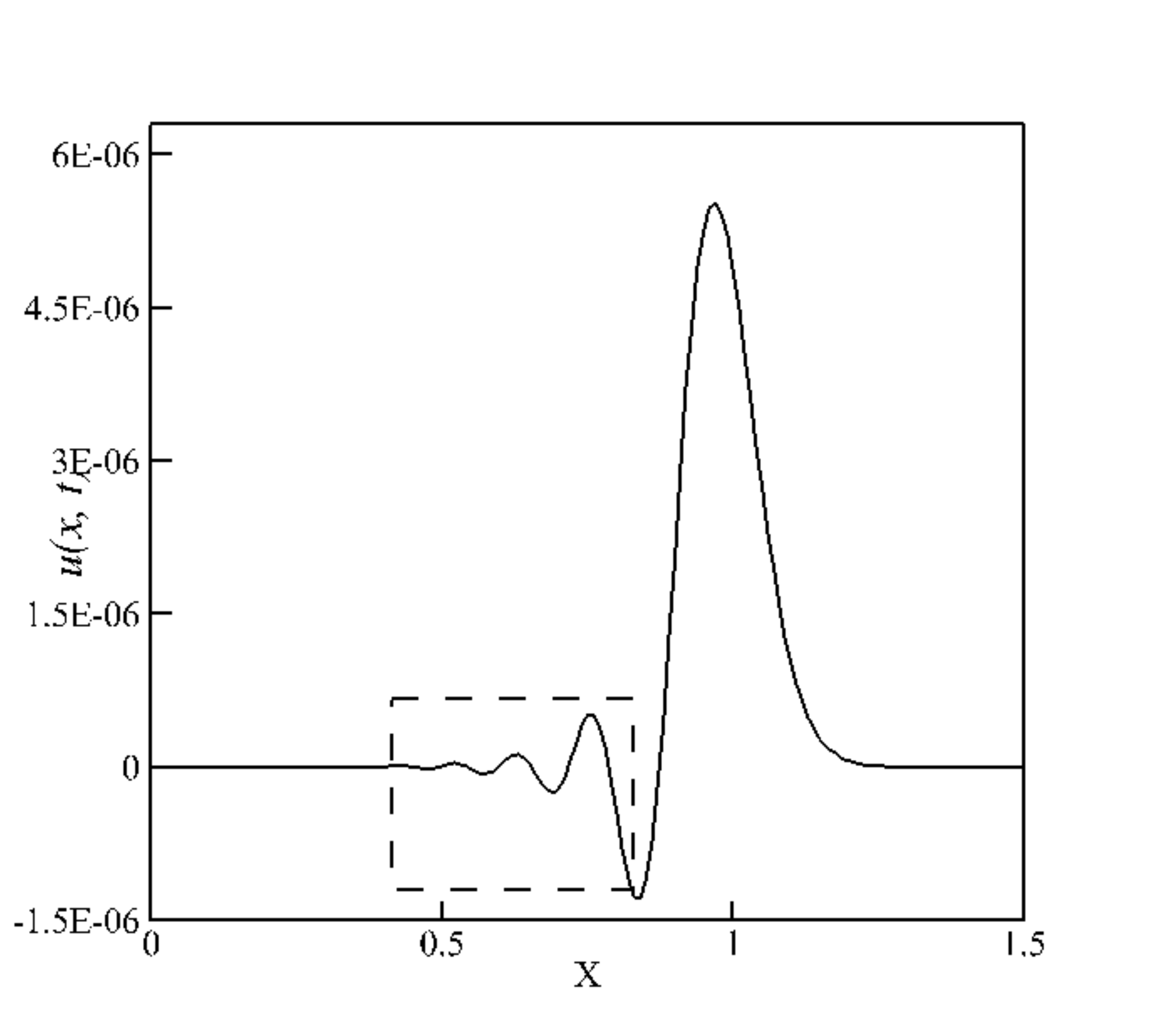}
 \vskip -5pt
 \caption{} \label{fig:Fig5d}
\end{subfigure}
\vskip -3.5pt
\begin{subfigure}{0.48\textwidth}
\includegraphics[width=0.96\linewidth, height=0.8\linewidth]{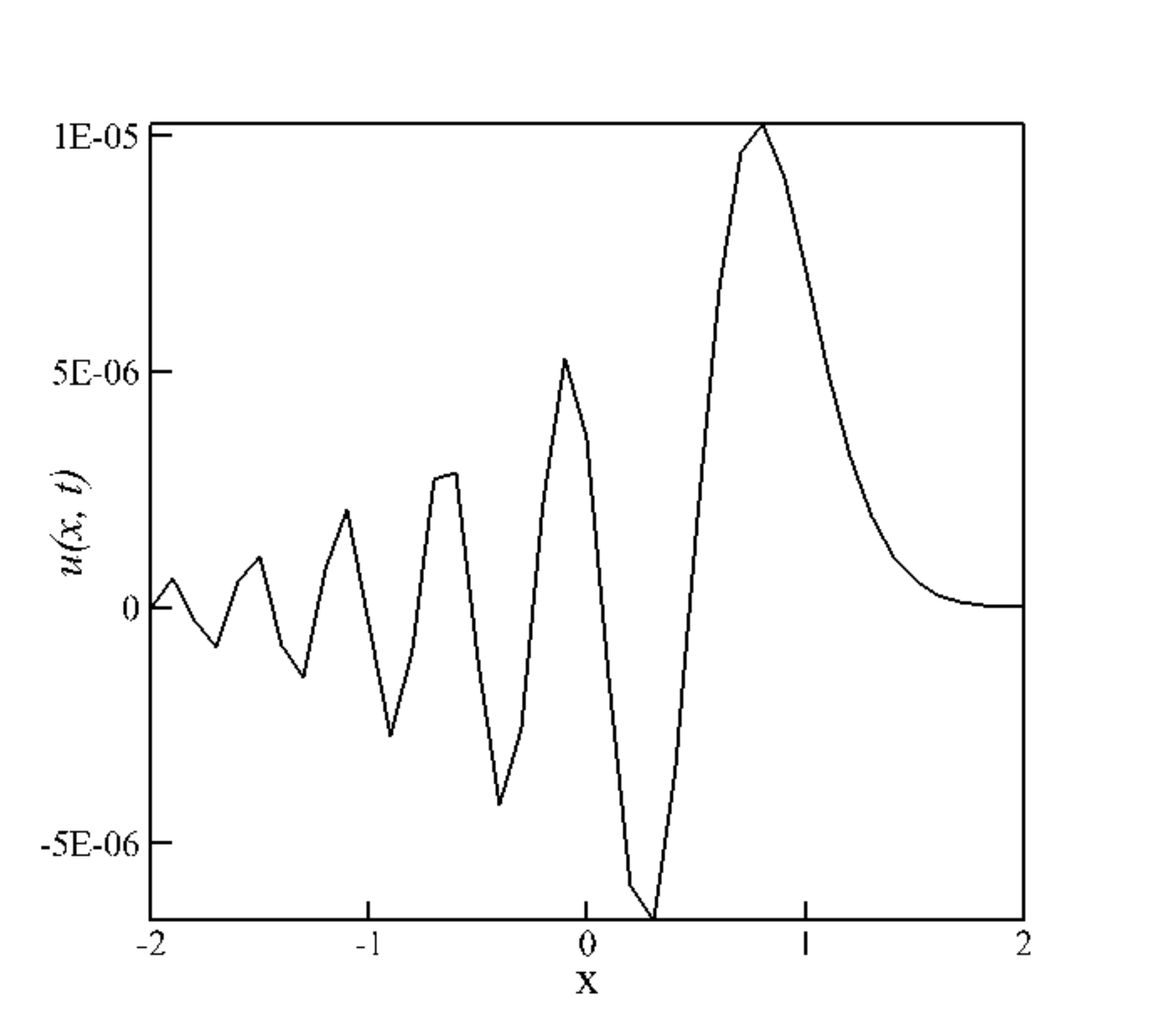}
\vskip -5pt
\caption{} \label{fig:Fig5e}
\end{subfigure}
\begin{subfigure}{0.48\textwidth}
 \includegraphics[width=0.96\linewidth, height=0.8\linewidth]{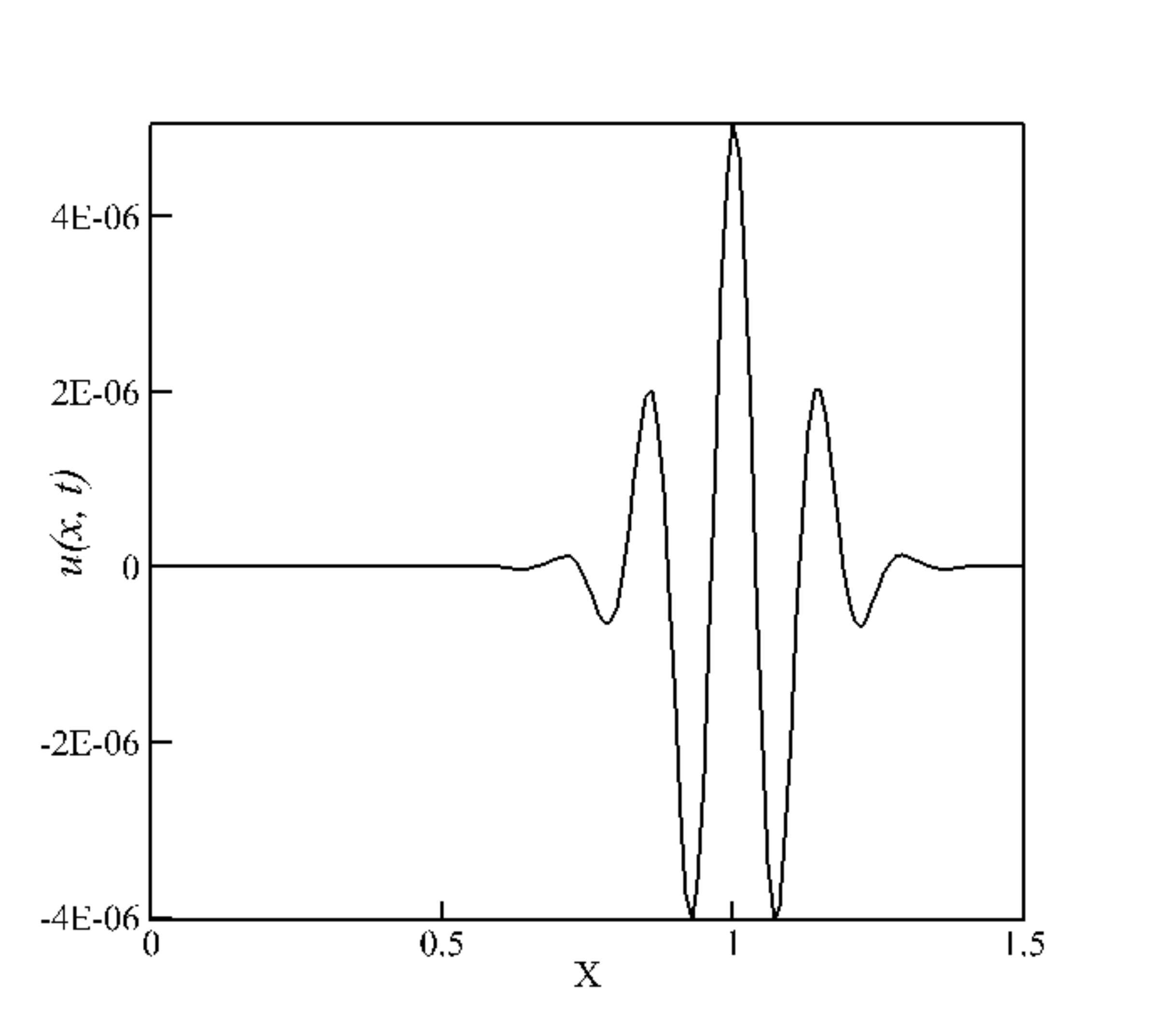}
 \vskip -5pt
 \caption{} \label{fig:Fig5f}
\end{subfigure}
\vskip -5pt
\caption{(a) Exact solution of the equation~\eqref{eqn:ADR} along with the initial condition~\eqref{eqn:IC-ADR}, (b) Fourier transform of the exact solution, (c) solution at $\gamma=50$ computed via the Explicit-OUCS3-CD2 scheme (\S \ref{subsec:scheme1}), (d) solution at $\gamma=10000$ computed using the Explicit-OUCS3-CD2 scheme where the boxed region highlights the presence of $q-$waves, (e) solution at $\gamma=10000$ computed by utilizing the Explicit-OUCS3-CD2 scheme on a coarser grid, $N=101$, and (f) solution at $\gamma=50$ computed by using the IMEX-NCCD scheme (\S \ref{subsec:scheme4}). The other parameters are fixed at $c=0.1,  \nu=10^{-4}, \lambda=-1.0, x_0=0.0, k_0h=0.5, L=5.0$ and $N=1001$ (except in figure~\ref{fig:Fig5e}). All numerical results are shown at non-dimensional simulation time $t=10.0$.}\label{fig:Fig5}
\end{figure}

Figure~\ref{fig:Fig5c}, \ref{fig:Fig5f} highlights the numerical solution of the 1D ADR equation~\eqref{eqn:ADR} for an initial wave-packet of width $\gamma=50$ (equation~\eqref{eqn:IC-ADR}) using the Explicit scheme (\S \ref{subsec:scheme1}) and the IMEX-NCCD scheme (\S \ref{subsec:scheme4}), respectively. At initial excitation $k_0h = 0.5$ (and $N_c = 0.1$), the Explicit scheme estimates the propagation characteristics of the signal poorly (i.~e., $V_g=0.7078$) while the IMEX-NCCD scheme predicts this value relatively accurately ($V_g=1.0013$). Another source of error predicted by the error propagation equation~\eqref{appendA:ErrorADR} is the absolute error in the phase speed, which is observed to be $c^{\text{err}}_{\text{phase}}=0.0956 (4.3686 \times 10^{-4}$) for the Explicit-OUCS3-CD$_2$ (IMEX-NCCD) scheme. A substantial phase error of the Explicit scheme is responsible for the dispersion effects resulting in an asymmetry in the computed solution of equation~\eqref{eqn:ADR} (figure~\ref{fig:Fig5c}) albeit the exact solution is symmetric (figure~\ref{fig:Fig5a}). The amplification ratio of the initial excitation is estimated to be $G=0.9896 / 0.9999$ for the Explicit-OUCS3-CD$_2$ / IMEX-NCCD scheme, respectively. We note that the solution obtained by IMEX-NCCD scheme has negligible $q-$waves (figure~\ref{fig:Fig5f}). The presence of $q-$waves implies that the energy of the $p-$waves (or physical waves) is siphoned off to the $q-$waves ensuing a diminished amplitude of the numerical solution (e.~g., compare the amplitude of the computed solution of the Explicit scheme (figure~\ref{fig:Fig5c}) versus the IMEX-NCCD scheme (figure~\ref{fig:Fig5f})). The fact that the amplification ratio of the IMEX-NCCD scheme remains very close to one indicates the near absence of $q-$waves. %in this case (refer figure~\ref{fig:Fig5f}).

We summarize our discussion by indicating two sources of error which are particularly perplexing in the DNS of ADR equations, the existence of $q-$waves for those set of numerical parameters for which the spatiotemporal discretization is neutrally stable. In such a situation the $q-$waves do not attenuate and have to be eliminated by deploying an explicit filter~\cite{TKS2009LES}. Another aspect of dispersion error is related to the Gibbs' phenomena which occurs as a consequence of sharp discontinuity in the solution and which causes fictitious oscillations~\cite{TKS2007,TKS2014}, a problem which can be remedied using the high accuracy DRP methods. Both these sources of error are discernible in nonlinear ADR equation whose solution (for a particular case) is explored next.

%While comparing the numerical methods in \S \ref{sec:spectral}, one notes that the limiting value of the wavenumber beyond which the $q-$waves are created is significantly higher in the IMEX-NCCD method ($kh > 2.37$, refer figure~\ref{fig:Fig4b}). 

%%%%%%%%%%%%%%%%%%%%%%%%%%%%%%%%%%%%%%%%%%%%%%%%
\section{IMEX method for chemotaxis model} \label{sec:PKS}
In this section we demonstrate the suitability of the IMEX-NCCD method developed in \S \ref{subsec:scheme4} for solving the PKS chemotaxis model~\cite{Patlak1953}. Chemotaxis refers to the mechanism by which cellular motion occurs in response to an external, chemical stimulus~\cite{Keller1970}. The PKS model is governed by a system of ADR equations, which in the classical 2D case reads as, 
\begin{align}
& \rho_t + \nabla \cdot (\chi \rho \nabla c) = \triangle \rho, \nonumber \\
& c_t = \triangle c - c + \rho, \qquad \qquad (x, y) \in \Omega \subset \mathbb{R}^2, \quad t > 0
\label{eqn:PKS}
\end{align}
where $\rho(x, y, t), c(x, y, t)$ denotes the cell density and the chemoattractant concentration, respectively. $\chi > 0$ is the chemotactic sensitivity parameter. $\nabla = (\frac{\partial }{\partial x}, \frac{\partial }{\partial y})$ and $\triangle = \frac{\partial^2}{\partial x^2} + \frac{\partial^2}{\partial y^2}$ represents the 2D gradient and the Laplacian operators, respectively. A common property of all chemotaxis systems is their ability to accurately capture concentrations that mathematically result in a rapid growth in small neighborhoods of the solution. For example, the solution of equation~\eqref{eqn:PKS} may `blow up' in finite time, provided that the initial mass of the cells, $\int_{\Omega} \rho(x, y, 0) dx dy$, is above a critical threshold which is $\nicefrac{8\pi}{\chi}$ for radially symmetric cases~\cite{Keller1971}. Steep gradients and spikes may give rise to nonphysical oscillations if the numerical scheme is not guaranteed to satisfy the discrete maximum principle (DMP), resulting in negative cell densities~\cite{Strehl2013}. To capture the singular, spiky behavior of the solution of~\eqref{eqn:PKS}, a high accuracy, positivity-preserving, Finite Volume (FV) scheme (satisfying DMP) is deployed by rewriting equation \eqref{eqn:PKS} as follows, %IMEX-NCCD scheme (satisfying DMP) is deployed by rewriting equation \eqref{eqn:PKS} as follows,
\begin{align}
& \rho_t + (\chi \rho u)_x + (\chi \rho v)_y = \triangle \rho, \nonumber \\
& c_t = \triangle c - c + \rho, \qquad \qquad \qquad \qquad \quad u = c_x, \,\, v = c_y,
\label{eqn:PKS-2D}
\end{align}
where $u$ and $v$ are the chemotactic velocities. The subscripts $x, y$ in equation~\eqref{eqn:PKS-2D} denote the partial derivative with the respective variables. We consider the model~\eqref{eqn:PKS-2D} in a square domain, $\Omega$, and introduce a Cartesian mesh consisting of cells $I_{i, j}=[x_{i-\frac{1}{2}}, x_{i+\frac{1}{2}}] \times [y_{j-\frac{1}{2}}, y_{j+\frac{1}{2}}]$ which are assumed to be of uniform size ($hk$), such that $x_{i+\frac{1}{2}} - x_{i-\frac{1}{2}} \equiv h$ for all $i$ and $y_{j+\frac{1}{2}} - y_{j-\frac{1}{2}} \equiv k$ or all $j$. On this mesh, a general FV scheme for the PKS system~\eqref{eqn:PKS-2D} will have the following form, %-IMEX-NCCD scheme for the PKS system~\eqref{eqn:PKS-2D} will have the following form,
\bseq
\begin{align}
\frac{d \rho_{i,j}}{d t} & = -\frac{\mathcal{P}_{i+\frac{1}{2},j}-\mathcal{P}_{i-\frac{1}{2},j}}{h} -\frac{\mathcal{Q}_{i,j+\frac{1}{2}}-\mathcal{Q}_{i,j-\frac{1}{2}}}{h} + \triangle_{i,j} \rho, \label{subeqn:rho} \\
\frac{d \rho_{i,j}}{d t} & = \triangle_{i,j} c - c_{i, j} + \rho_{i, j}. \label{subeqn:c}
\end{align}
\label{eqn:discretePKS}
\eseq
The values $\rho_{i,j}$ and $c_{i,j}$ are the approximate point variable values at the cell center, $\triangle_{i,j}$ is the discrete Laplacian operator, $\mathcal{P}_{i+\frac{1}{2}, j} = \chi \rho_{i+\frac{1}{2}, j} u_{i+\frac{1}{2}, j}$ and $\mathcal{Q}_{i, j+\frac{1}{2}} = \chi \rho_{i, j+\frac{1}{2}} v_{i, j+\frac{1}{2}}$ are the numerical fluxes in the x- and the y- direction, respectively, and these are evaluated at the cell edges. The derivatives, $u_{i+\frac{1}{2},j}$ and $v_{i,j+\frac{1}{2}}$ are approximated by appropriate compact / explicit finite difference (FD) scheme at the cell centers and then linearly extrapolated at the cell edges (refer \S \ref{sec:spectral} for the details of the spatial discretization utilized). The point values $\rho_{i+\frac{1}{2},j}$ and $\rho_{i,j+\frac{1}{2}}$ are computed as follows,
%
%spatially implicit NCCD scheme at the cell centers and then extrapolated at the cell edges, whereas the point values $\rho_{i+\frac{1}{2},j}$ and $\rho_{i,j+\frac{1}{2}}$ are computed as follows,
%
\vskip -10pt
\begin{minipage}{0.5\textwidth}
\beq
\rho_{i+\frac{1}{2},j}=\begin{cases}
\rho_{i,j} + \frac{h}{2} (\rho_x)_{i,j} \quad \qquad \text{if}\,\,u_{i+\frac{1}{2},j} > 0, \\
\rho_{i+1,j} - \frac{h}{2} (\rho_x)_{i+1,j} \quad \text{if}\,\, u_{i+\frac{1}{2},j} \le 0, \nonumber
\end{cases}
\eeq
\end{minipage}
\begin{minipage}{0.5\textwidth}
\beq
\rho_{i,j+\frac{1}{2}}=\begin{cases}
\rho_{i,j} + \frac{k}{2} (\rho_y)_{i,j} \quad \qquad \text{if}\,\, v_{i,j+\frac{1}{2}} > 0, \\
\rho_{i,j+1} - \frac{k}{2} (\rho_y)_{i,j+1} \quad \text{if}\,\, v_{i,j+\frac{1}{2}} \le 0,
\end{cases}
\label{eqn:discreteRho}
\eeq
\end{minipage}
In order to ensure that the point values in equation~\eqref{eqn:discreteRho} are nonnegative, the slopes are calculated adaptively using
\bseq
\begin{align}
(\rho_x)_{i,j} &=\begin{cases}
(\rho_x)_{i,j} \,\,\,\,\,\,\,\quad\quad\qquad\qquad\qquad\qquad\qquad\qquad \text{if}\,\, \rho_{i,j} \pm \frac{h}{2} (\rho_x)_{i,j} \ge 0, \\
\text{minmod}(\theta (\rho_x)_{i+\frac{1}{2},j}, \,\, (\rho_x)_{i,j}, \,\, \theta (\rho_x)_{i-\frac{1}{2},j}) \quad \text{otherwise}, %\nonumber
\end{cases} \\
(\rho_y)_{i,j} &=\begin{cases}
(\rho_y)_{i,j} \,\,\,\,\,\,\,\quad\quad\qquad\qquad\qquad\qquad\qquad\qquad \text{if}\,\, \rho_{i,j} \pm \frac{k}{2} (\rho_y)_{i,j} \ge 0, \\
\text{minmod}(\theta (\rho_y)_{i,j+\frac{1}{2}}, \,\, (\rho_y)_{i,j}, \,\, \theta (\rho_y)_{i,j-\frac{1}{2}}) \quad \text{otherwise}, %\nonumber
\end{cases}
\end{align}
\label{eqn:PPS}
\eseq
where the positivity-preserving generalized minmod limiter, given by
\beq
\text{minmod}(r_1, r_2, \ldots) = \begin{cases}
\text{min}(r_1, r_2, \ldots) \quad \text{if} \quad r_i > 0 \quad \forall i, \\
\text{max}(r_1, r_2, \ldots) \quad \text{if} \quad r_i < 0 \quad \forall i, \\
0 \,\,\,\qquad\qquad\qquad \text{otherwise},
\end{cases}
\eeq
ensures the positivity of the reconstructed point values. The parameter $\theta$ in equation~\eqref{eqn:PPS} controls the amount of numerical viscosity present in the scheme, larger values of $\theta$ correspond to less dissipative but usually more oscillatory reconstructions. An explicit FV-FD numerical scheme for solving the 2-D Keller-Segel model is outlined in~\cite{Chertock2017}.
%spatially explicit version of the above mentioned numerical scheme, for solving the 2-D Keller-Segel model is outlined in~\cite{Chertock2019}.

%%%%%%%%%%%%%%%%%%%%%%%%%%%%%%%%%%%%%%%%%%%%%%%%
\subsection{Numerical results} \label{subsec:NR}
\begin{figure}[htbp]
\centering
\begin{subfigure}{0.48\textwidth}
\includegraphics[width=0.96\linewidth, height=0.8\linewidth]{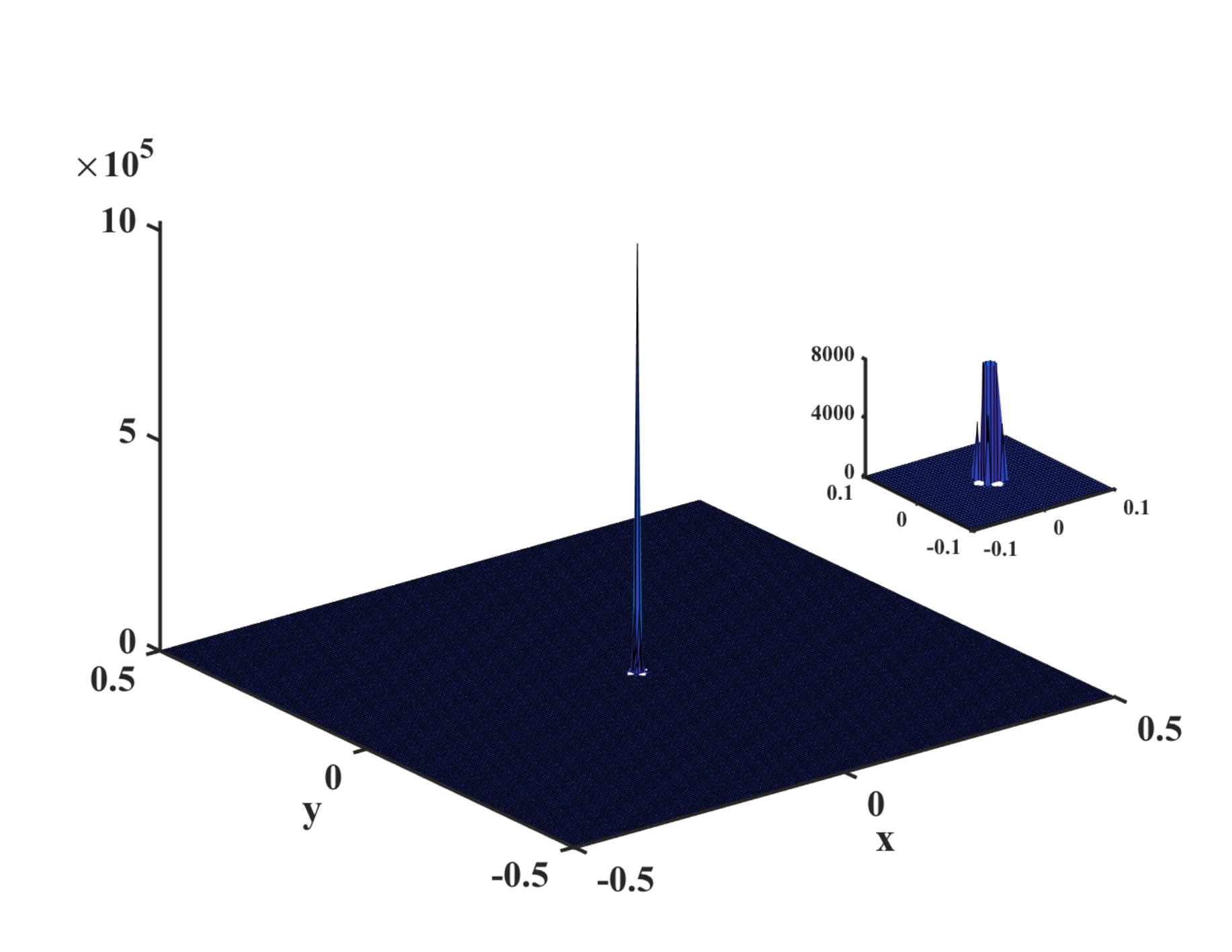}
\vskip -5pt
\caption{} \label{fig:Fig6a}
\end{subfigure}
\begin{subfigure}{0.48\textwidth}
\includegraphics[width=0.96\linewidth, height=0.8\linewidth]{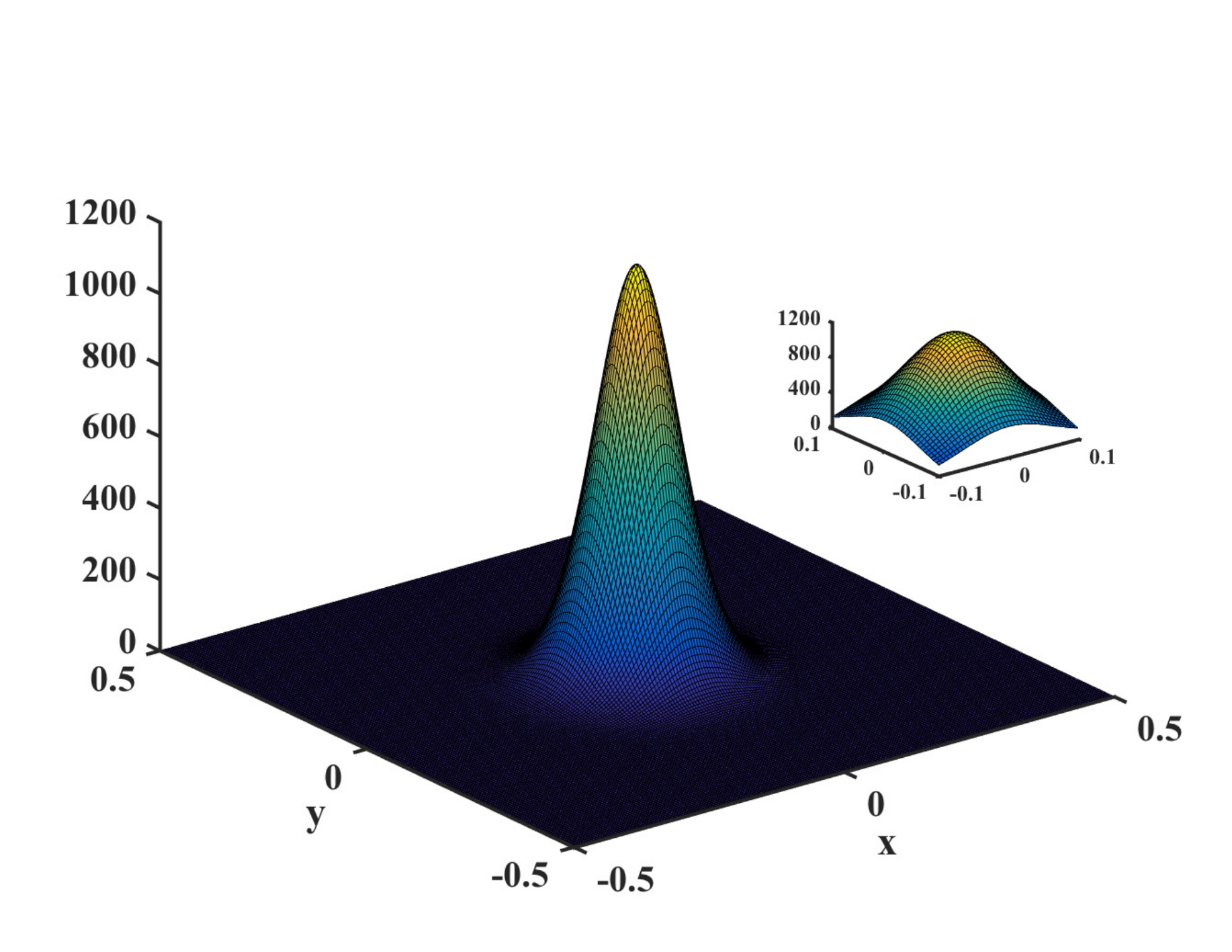}
\caption{} \label{fig:Fig6b}
\end{subfigure}
\begin{subfigure}{0.48\textwidth}
 \includegraphics[width=0.96\linewidth, height=0.8\linewidth]{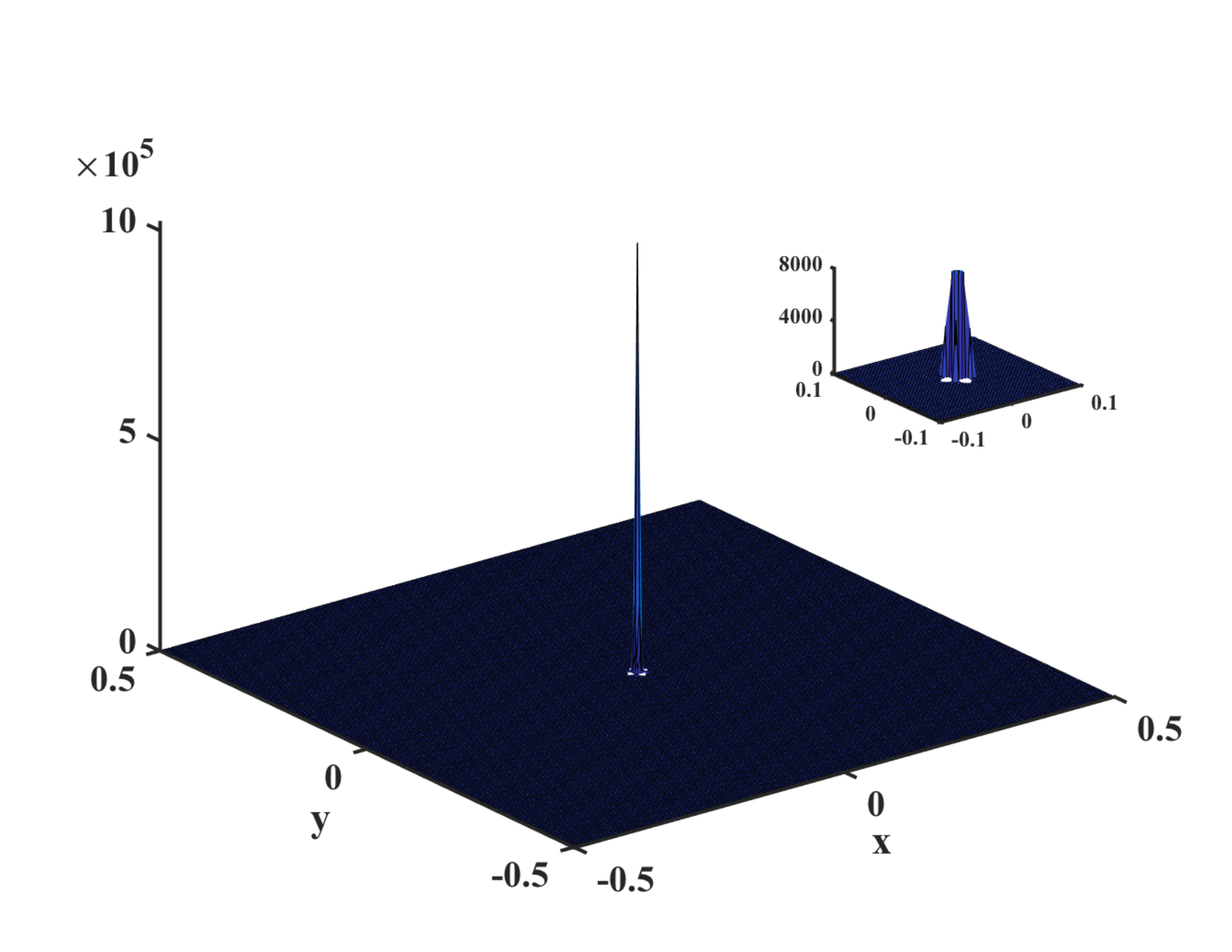}
 \vskip -5pt
 \caption{} \label{fig:Fig6c}
\end{subfigure}
\begin{subfigure}{0.48\textwidth}
\includegraphics[width=0.96\linewidth, height=0.8\linewidth]{Fig6c.pdf}
\caption{} \label{fig:Fig6d}
\end{subfigure}
\caption{The cell density, $\rho$ (equation~\eqref{eqn:PKS}), computed by (a) FV-Explicit-OUCS3-CD$_2$ method at $T = 10^{-5}$, (b) FV-IMEX-NCCD method at $T = 10^{-5}$, (c) FV-Explicit-OUCS3-CD$_2$ method at $T = 5\times10^{-3}$ and (d) FV-IMEX-NCCD method at $T = 5\times10^{-3}$ on a uniform mesh with $\triangle x = \triangle y = 1/200$. }\label{fig:Fig6}
\end{figure}

An initial exploration to test the efficacy of the IMEX-NCCD method was conducted using an example taken from~\cite{Chertock2017}. We consider the initial-boundary value problem for the PKS system~\eqref{eqn:PKS} with $\chi=30.0$ and $\theta=1.0$ (refer \eqref{eqn:PPS}) on a square domain $[-\frac{1}{2}, \frac{1}{2}] \times [-\frac{1}{2}, \frac{1}{2}]$ subject to a radially symmetric gaussian initial condition:
\beq 
\rho(x, y, 0) = 1000 e^{-100(x^2 +y^2)}, \qquad c(x, y, 0) = 500 e^{-50(x^2 +y^2)}.
\eeq
As it was highlighted in~\cite{Chertock2017}, the solution with the above mentioned initial data is expected to develop a $\delta-$like singularity at the center of the computational domain in a very short time. Figure~\ref{fig:Fig6} compares the cell densities computed via the FV-Explicit-OUCS3-CD$_2$ method (figures~\ref{fig:Fig6a}, \ref{fig:Fig6c}) versus the FV-IMEX-NCCD scheme (figures~\ref{fig:Fig6b}, \ref{fig:Fig6d}) on a uniform mesh with $h = k = 1/200$, at two different computational times $T = 10^{-5}$ (Figures~\ref{fig:Fig6a}, \ref{fig:Fig6b}) and $T = 5 \times 10^{-3}$ (Figures~\ref{fig:Fig6c}, \ref{fig:Fig6d}), respectively. In both of these cases the zero Neumann boundary conditions are used, which is implemented using the standard ghost point technique when the CD$_2$ discretization is used. Further, the time steps $dt=10^{-8}$ and $dt=10^{-6}$ were chosen for the respective methods to ensure the positivity of the computed densities. As outlined in figure~\ref{fig:Fig6} although both methods capture the `needle-like' structure, the blowup phenomena is better resolved by the IMEX-NCCD method since the solution obtained via this method is oscillation-free (compare the figure insets of figure~\ref{fig:Fig6a} versus figure~\ref{fig:Fig6b} or that of figure~\ref{fig:Fig6c} versus figure~\ref{fig:Fig6d}). This is because, centered approximations (e.~g.~the CD$-2$ discretization) for hyperbolic PDEs lead to well-known difficulties with nonphysical oscillations (or Gibbs' phenomena) arising near discontinuities or steep gradients~\cite{LeVeque1992}. A detailed study relating the origin of these oscillations with the time-step of the numerical methods will be reported later.

%%%%%%%%%%%%%%%%%%%%%%%%%%%%%%%%%%%%%%%%%%%%%%%%
\section{Summary and conclusion} \label{sec:conclusion}
This introductory investigation compares and reports the dispersion error of four, two-time level, DRP schemes using the novel error propagation equation for 1D linear ADR equation (details given in \S \ref{sec:appendA}). \S \ref{sec:spectral} presented a comprehensive spectral analysis of these four numerical methods, highlighting regions of temporal stability, positive scaled group velocity and minimal phase errors in the wavenumber ($kh$) - CFL number ($N_c$) plane. The role of the initial conditions, inadequate grid resolution or the choice of the numerical methods in the generation of the dispersion errors via $q-$waves is outlined in \S \ref{sec:ED}. The superior spectral resolution of the IMEX-NCCD method is then benchmarked via the solution of the nonlinear, parabolic-elliptic PKS chemotaxis model in \S \ref{sec:PKS}. The blowup phenomena in the numerical solution of the PKS model was accurately captured by the IMEX-NCCD scheme, devoid of any amplitude fluctuations. In a nutshell, the numerical solution obtained by the IMEX-NCCD method has negligible dispersion error due to $q-$waves and is ideally suited to arrest the oscillations via the Gibbs' phenomena.

%%%%%%%%%%%%%%%%%%%%%%%%%%%%%%%%%%%%%%%%%%%%%%%%
\section{Appendix A: Derivation of the error propagation equation for 1D linear ADR equation} \label{sec:appendA}
The analysis of the space-time discretization of the 1D linear ADR equation~\eqref{eqn:ADR}, with initial condition given by, $u(x, 0) = f(x)$ for $x \in (-\infty, \infty)$, is used as a model for PDE replicating multiscale processes. Using the spectral representation $u(x, t) = \int U(k, t) e^{{\it i}k x} d k$, equation~\eqref{eqn:ADR} is transformed as
\beq
U_t + (\nu k^2 + {\it i}k c - \lambda) U = 0,
\label{appendA:spectralODE}
\eeq
where $k \in (-\infty, \infty)$ is the wavenumber. Using the transformed initial condition, $f(x) = \int A_0(k) e^{{\it i}k x} d k$, one obtains the solution of equation~\eqref{appendA:spectralODE} as
\beq
U(k, t) = A_0(k) e^{-(\nu k^2 + {\it i}k c - \lambda)t}.
\label{appendA:U}
\eeq
If we represent $u(x, t)$ by its Fourier-Laplace transform as, $u(x, t) = \int \int \hat{U}(k, \omega) e^{{\it i}(k x - \omega t)}d k d \omega$, then the dispersion relation for equation~\eqref{eqn:ADR} in the physical plane is
\beq
\omega = k c - {\it i} \nu k^2 + {\it i} \lambda,
\label{appendA:DRPphy}
\eeq
and using equation~\eqref{appendA:U}, one obtains the time-amplification factor of the PDE, 
\beq
G_{\text{exact}} = \nicefrac{U(k, t+\triangle t)}{U(k,t)} = e^{-(\nu k^2 + {\it i}k c - \lambda)\triangle t}. 
\label{appendA:Gexact}
\eeq
Errors arise due to discrete computing of equation~\eqref{eqn:ADR} by numerical methods, as described next. 

Suppose the unknown, $u^N(x_j, t)$, at the $j$th node of a uniformly spaced grid of spacing $h$ can be represented in the spectral plane by $u^N(x_j, t) = \int U(k, t) e^{{\it i}k x_j} d k$ and the initial condition by, $u^N(x_j, t = 0) = u_0^N = \int A_0(k) e^{{\it i}k x_j} d k$. With $k$ as the independent variable, it is concluded that the coefficients, $c, \nu, \lambda$, become $k-$dependent numerically~\cite{TKS2007}, culminating in the numerical dispersion relation for equation~\eqref{eqn:ADR} as
\beq
\omega^N = k c^N - {\it i} \nu^N k^2 + {\it i} \lambda^N,
\label{appendA:DRPnum}
\eeq
where the superscript `N' denotes the corresponding (complex) numerical values. The general numerical solution at any time instant is given by
\beq
u^N_j = \bigintsss A_0(k) [G_{\text{num}}]^{\frac{t}{\triangle t}} e^{{\it i}k x_j} d k,
\label{appendA:Unum}
\eeq
where the numerical amplification factor, $G_{\text{num}}$, is obtaining by substituting the respective numerical values of the coefficients, $c^N, \nu^N, \lambda^N$ in equation~\eqref{appendA:Gexact}. If we define error as $e(x, t) = u(x, t) - u^N$, then from equation~\eqref{eqn:ADR}, we obtain the dynamics of error propagation as
\begin{align}
\frac{\partial e}{\partial t} + c \frac{\partial e}{\partial x} - \nu \frac{\partial^2 e}{\partial x^2} - \lambda e &= -\frac{\partial u^N}{\partial t} -c\frac{\partial u^N}{\partial x} + \nu \frac{\partial^2 u^N}{\partial x^2} + \lambda u^N \nonumber \\
&= \bigintsss {\it i} k A_0(k) c_{\text{exact}} \left[ 1 - \frac{c_{\text{num}}}{c_{\text{exact}}}\right] [G_{\text{num}}]^{t/\triangle t} e^{{\it i}k x} dk %\nonumber \\
%
%&= \bigintsss \!\!\!\!A_0(k) ({\it i} k c + \nu k^2 \!\!-\!\! \lambda) \!\!\!\left[ 1 \!\!+\!\! \frac{h}{c \triangle t}\!\!\left( \frac{\ln |G_{\text{num}}| \!-\! {\it i} \beta}{\frac{\nu (k h)^2 - \lambda h^2}{c h} + {\it i}k h} \right) \!\!\right] \!\![G_{\text{num}}]^{t/\triangle t} e^{{\it i}k x} dk,
\label{appendA:ErrorADR}
\end{align}
where the numerical phase speed, $c_{\text{num}}=\nicefrac{\omega^N}{k}$ (equation~\eqref{appendA:DRPnum}), and the physical phase speed, $c_{\text{exact}}=\nicefrac{\omega}{k}$ (equation~\eqref{appendA:DRPphy}), are given by
\begin{align}
\frac{c_{\text{num}}}{c_{\text{exact}}} &= -\frac{h}{c \triangle t}\!\!\left( \frac{\ln |G_{\text{num}}| \!-\! {\it i} \beta}{\frac{\nu (k h)^2 - \lambda h^2}{c h} + {\it i}k h} \right) \nonumber \\
c_{\text{exact}} &= c - {\it i} \nu k - \nicefrac{\lambda}{{\it i} k},
\label{appendA:phase}
\end{align}
and the numerical phase shift is given by $\beta=-\tan^{-1}\!\!\left[\frac{(G_{\text{num}})_{\text{imag}}}{(G_{\text{num}})_{\text{real}}}\right]$. Note that the right hand side of equation~\eqref{appendA:ErrorADR} is zero (which is the case of von Neumann stability analysis~\cite{Charney1950}) provided the error in the phase speed, $1 - \nicefrac{c_{\text{num}}}{c_{\text{exact}}} = 0$ and/or the numerical amplification factor, $G_{\text{num}} = G_{\text{exact}}$ (refer equation~\eqref{appendA:Gexact}). Finally the group velocity, V$_g$, can be calculated from the dispersion relation by V$_g = \frac{\partial \omega}{\partial k}$ and on further simplification yields the ratio,
\beq
\frac{V_{g,\text{num}}}{V_{g,\text{exact}}} = \frac{h}{c \triangle t} \frac{d \beta}{d (kh)},
\label{appendA:Vg}
\eeq
where $V_{g,\text{num}}$, $V_{g,\text{exact}}$ are the numerical and the physical group velocity, respectively.

\paragraph{Acknowledgements:} The first author acknowledges the financial support of the grant DST ECR / 2017 / 000632 which helped him to run the numerical simulations. The second author acknowledges the financial support of the grant DST MTR / 2017 / 000017.

\bibliographystyle{siam}
\bibliography{JSC}

\end{document}